\documentclass[11pt,reqno]{amsart}
\usepackage[T1]{fontenc}
\usepackage[letterpaper,margin=1.25in]{geometry}
\usepackage{verbatim}
\usepackage{amssymb,amsmath,amsthm}
\usepackage{booktabs}
\usepackage{float}
\usepackage[colorlinks=true,linkcolor=blue,citecolor=blue,urlcolor=blue]{hyperref}
\usepackage{cleveref}
\DeclareMathOperator{\rk}{rk}
\DeclareMathOperator{\Aut}{Aut}
\DeclareMathOperator{\End}{End}
\DeclareMathOperator{\Gal}{Gal}
\DeclareMathOperator{\Cl}{Cl}
\DeclareMathOperator{\gon}{gon}
\DeclareMathOperator{\ddiv}{div}

\DeclareMathOperator{\PSL}{PSL}
\DeclareMathOperator{\SL}{SL}
\DeclareMathOperator{\ord}{ord}
\DeclareMathOperator{\Frob}{Frob}
\DeclareMathOperator{\Ver}{Ver}

\newcommand{\Q}{\mathbb{Q}}
\newcommand{\Z}{\mathbb{Z}}
\newcommand{\C}{\mathbb{C}}
\newcommand{\F}{\mathbb{F}}
\newcommand{\PP}{\mathbb{P}}
\newcommand{\Qbar}{\overline{\Q}}
\newcommand{\tors}{\mathrm{tors}}
\newcommand{\torz}[1]{\Z/#1\Z}
\newcommand{\tg}[2]{\Z/#1\Z\times\Z/#2\Z}
\newcommand{\dm}[1]{\langle #1\rangle}
\newcommand{\ClCusp}{\Cl^{c}}
\newcommand{\Cart}{\mathcal{C}}

\newtheorem{theorem}{Theorem}[section]
\newtheorem{lemma}[theorem]{Lemma}
\newtheorem{proposition}[theorem]{Proposition}

\theoremstyle{definition}
\newtheorem{remark}[theorem]{Remark}

\newcommand{\lmfdbec}[3]{\href{https://www.lmfdb.org/EllipticCurve/Q/#1/#2/#3}{#1.#2#3}}
\newcommand{\lmfdcharorbit}[2]{\href{https://www.lmfdb.org/Character/Dirichlet/#1/#2}{#1.#2}}
\newcommand{\lmfdbnewform}[4]{\href{https://www.lmfdb.org/ModularForm/GL2/Q/holomorphic/#1/#2/#3/#4}{#1.#2.#3.#4}}
\newcommand{\lmfdbnf}[1]{\href{https://www.lmfdb.org/NumberField/#1}{#1}}
\newcommand{\xpoint}[1]{\href{https://x1mn-sporadic.github.io/point.html?id=#1}{#1}}
\newcommand{\githubbare}[1]{\href{https://github.com/F-Najman/sextic_torsion/blob/dbcd360dbe576dadfd0da581d09d745e40429a45/#1}{\path{#1}}}

\newcommand{\gitlink}[2]{\href{https://github.com/F-Najman/sextic_torsion/blob/dbcd360dbe576dadfd0da581d09d745e40429a45/#1}{#2}}
\newcommand{\githubdir}[1]{\href{https://github.com/F-Najman/sextic_torsion/tree/dbcd360dbe576dadfd0da581d09d745e40429a45/#1}{\path{#1}}}

\title[Torsion of elliptic curves over sextic fields]
      {Classification of torsion of elliptic curves\\ over sextic fields}
\author{Filip Najman}
\address{University of Zagreb Faculty of Science, Department of Mathematics,
Bijeni\v cka cesta 30, 10000 Zagreb, Croatia}
\email{fnajman@math.hr}
\thanks{The author was supported by the Croatian Science Foundation under the project no.\ IP-2022-10-5008, by the project ``Implementation of cutting-edge research and its application as part of the Scientific Center of Excellence for Quantum and Complex Systems, and Representations of Lie Algebras'', PK.1.1.10.0004, co-financed by the European Union through the European Regional Development Fund -- Competitiveness and Cohesion Programme 2021--2027, and by the European Union -- NextGenerationEU through the National Recovery and Resilience Plan 2021--2026, via an institutional grant of the University of Zagreb Faculty of Science, IK IA 1.1.3, Impact4Math.}
\subjclass[2020]{11G05, 11G18, 14G05}
\keywords{elliptic curves, torsion, sextic fields, modular curves, sporadic points}

\begin{document}

\begin{abstract}
We determine all the groups that appear as the torsion group of an elliptic curve over a number field of degree $6$. Apart from the $44$ groups that occur infinitely often, which were determined by Derickx and Sutherland, exactly four groups occur: $\torz{25}$, $\torz{37}$, $\tg{4}{12}$ and $\tg{2}{24}$. The first three were known to occur, while the group $\tg{2}{24}$ is new. We also show that, up to isomorphism, each of the groups $\torz{37}$ and $\tg{2}{24}$ is the torsion group of a single elliptic curve over a single sextic field. The methods build on those developed by Derickx and Najman for quartic fields and by the author for quintic fields. Among the new ingredients are a formal immersion criterion at the prime $2$ which uses the Cartier operator, and arguments with several diamond operators and with cuspidal divisor classes on quotients of modular curves.
\end{abstract}

\date{\today}
\maketitle

\section{Introduction}

Let $E$ be an elliptic curve over a number field $K$. By the Mordell--Weil theorem
$E(K)\simeq E(K)_\tors\times \Z^r$, and by Merel's theorem \cite{merel} the torsion subgroup
$E(K)_\tors$ is bounded in terms of $d=[K:\Q]$ alone. Hence, for every $d\geq 1$ the set
\[
\Phi(d):=\{\,E(K)_\tors \;:\; [K:\Q]=d,\ E/K \text{ an elliptic curve}\,\}
\]
is finite. The set $\Phi(1)$ was determined by Mazur \cite{mazur}, building on work of
Kubert \cite{kubert}, and $\Phi(2)$ by Kamienny \cite{kamienny}, building on work of Kenku and
Momose \cite{KM}. Derickx, Etropolski, van Hoeij, Morrow and Zureick-Brown \cite{DEHMZ}
determined $\Phi(3)$, using earlier work of Jeon, Kim and Schweizer \cite{JKS}, Parent
\cite{parent} and Najman \cite{najman}. Derickx and Najman \cite{DN} determined $\Phi(4)$, using
earlier work of Jeon, Kim and Park \cite{JKP} and Bruin and Najman \cite{bruin_najman2016}.
Finally, $\Phi(5)$ was determined by the author \cite{N5}. In this paper we determine $\Phi(6)$.

Write $\Phi^\infty(d)\subseteq\Phi(d)$ for the subset of groups that occur as $E(K)_\tors$ for
infinitely many $\overline\Q$-isomorphism classes of elliptic curves $E$ over number fields $K$
of degree $d$. A group in $\Phi(d)\setminus\Phi^\infty(d)$ will be called a
\emph{sporadic torsion group in degree $d$}. This should not be confused with a \emph{sporadic
point} on a curve $X$, which is a closed point of some degree $e$ such that $X$ has only
finitely many closed points of degree at most $e$. A list of known sporadic and isolated points on
the curves $X_1(m,n)$ is maintained on the webpage \cite{x1mn}. We have $\Phi^\infty(d)=\Phi(d)$
for $d=1,2,4$.
The only sporadic torsion group in degree $3$ is $\torz{21}$ \cite{najman,DEHMZ}, and the
sporadic torsion groups in degree $5$ are $\torz{28}$, $\torz{30}$ and $\tg{2}{18}$
\cite{N5}. The set $\Phi^\infty(6)$ was determined by Derickx and Sutherland
\cite[Theorem 1.1]{DS}, the cyclic part being due to Derickx and van Hoeij
\cite[Theorem 3]{DvH}:
\begin{equation}\label{eq:phiinf6}
\begin{aligned}
\Phi^\infty(6)=\ &\{\torz{n}\;:\;1\leq n\leq 30,\ n\neq 23,25,29\}\\
              &\cup\ \{\tg{2}{2n}\;:\;1\leq n\leq 10\}\ \cup\ \{\tg{3}{3n}\;:\;1\leq n\leq 4\}\\
              &\cup\ \{\tg{4}{4},\ \tg{4}{8},\ \tg{6}{6}\}.
\end{aligned}
\end{equation}
This set has $44$ elements. Some sporadic torsion groups in degree $6$ were already known.
Van Hoeij \cite{vanHoeij} found an elliptic curve over a sextic field with a point of order
$25$ and one with a point of order $37$, and Derickx, Kamienny, Stein and Stoll proved that
every elliptic curve over a sextic field with a point of order $37$ has $j$-invariant $-9317$
\cite[Proposition 1.4]{DKSS}. Gonz\'alez-Jim\'enez and Najman \cite{GN} found that the base
change of the rational elliptic curve \lmfdbec{162}{d}{2} to a sextic field has torsion group
$\tg{4}{12}$ (see also \cite[Remark 1]{DGJ}).

Our main result is the following.

\begin{theorem}\label{thm:main}
Let $K$ be a number field of degree $6$ and let $E$ be an elliptic curve over $K$. Then
$E(K)_\tors$ is isomorphic to one of the following $48$ groups:
\begin{alignat*}{2}
&\torz{n}, \qquad && n=1,\dots,22,\ 24,\dots,28,\ 30,\ 37,\\
&\tg{2}{2n}, \qquad && n=1,\dots, 10,\ 12,\\
&\tg{3}{3n}, \qquad && n=1,\dots, 4,\\
&\tg{4}{4n}, \qquad && n=1,2,3,\\
&\tg{6}{6}. &&
\end{alignat*}
Equivalently, $\Phi(6)=\Phi^\infty(6)\cup\{\torz{25},\torz{37},\tg{4}{12},\tg{2}{24}\}$.
\end{theorem}

The four sporadic groups are described by the following theorem, which is proved in
\Cref{sec:sporadic}. For two of them we determine all the elliptic curves in question.

\begin{theorem}\label{thm:sporadic}
Let $K$ be a number field of degree $6$ and let $E/K$ be an elliptic curve.
\begin{enumerate}
\item If $E(K)$ has a point of order $37$, then $E(K)_\tors\simeq\torz{37}$, and up to
      isomorphism of the pair $(K,E)$, the field $K$ is the cyclic sextic field $K_{37}$ and $E$ is the curve $E_{37}$ of \Cref{prop:37}, a quadratic twist of \lmfdbec{1225}{b}{2}.
\item If $\tg{2}{24}\subseteq E(K)$, then $E(K)_\tors\simeq\tg{2}{24}$, and up to isomorphism
      of the pair $(K,E)$, the field $K$ is the field $K_{2,24}$ of discriminant $231852$ and $E$
      is the curve $E_{2,24}$ of \Cref{prop:2_24_points}.
\item There exist elliptic curves over sextic fields with torsion groups $\torz{25}$ and
      $\tg{4}{12}$, which are described in \Cref{prop:examples}.
\end{enumerate}
\end{theorem}

The curve in (2) is new. We do not determine all the elliptic curves over sextic fields with
torsion $\torz{25}$ or $\tg{4}{12}$. This is not needed for \Cref{thm:main}, since all the
groups properly containing $\torz{25}$ or $\tg{4}{12}$ can be ruled out directly.

\subsection*{Outline of the proof and of the paper}
As in \cite{DN} and \cite{N5}, the classification reduces to showing that certain modular
curves $X_1(m,n)$ have no non-cuspidal points of degree $6$, or more precisely no points
defined over a sextic field. In \Cref{sec:groups} we show that $66$ groups have to be
analyzed. Nine of them contain $\torz{37}$ and are ruled out by the explicit description of the
elliptic curves with a point of order $37$, and five more are reduced, by a $2$-isogeny, to
groups of the form $\tg{2}{2n}$. The others are treated with the following methods.

The global method of \cite[Proposition 4.3]{DN}, in the form used in \cite{N5}, eliminates
eleven groups $\torz{n}$ with large $n$ (\Cref{sec:global}). For six other large $n$ we use a
variant of the diamond argument of \cite[Proposition 7.2]{N5}, in which a single operator
$\dm a-1$ is replaced by a more general element of the group ring of the diamond operators.

When the Jacobian has rank $0$ we use the criterion of \cite[Proposition 5.2]{DN}, which
shows that the elliptic curve must have good reduction at some prime above a suitable $p$
(\Cref{sec:rank0}). When $p=2$ this often requires a large gonality, and we replace this
requirement by a condition on the dimensions of the Riemann--Roch spaces $L(2C)$ of the cuspidal
divisors $C$ on $X_{\F_2}$, which we check by computation. The cases in which the elliptic curve does have
good reduction at some prime above $p$ are eliminated by the Hecke sieve of
\cite[Proposition 5.5]{DN}, which we use both over $\Q$ and over $\Q(\zeta_m)$ for
$m=3,4,6$ (\Cref{sec:sieve}).

Several cases are settled by computations in the group generated by the cuspidal divisor
classes over a finite field (\Cref{sec:classes}). This determines all the sextic points on
$X_1(2,24)$, which is how we find the elliptic curves with torsion $\tg{2}{24}$ and rule out
$\torz{48}$. The groups $\torz{36}$, $\tg{2}{26}$, $\tg{3}{18}$, $\tg{2}{30}$ and $\torz{54}$
are also treated in this way, the last two by pushing divisors forward to quotient curves. For
$\torz{54}$ this leaves only elliptic curves whose reductions modulo $5$ have $j$-invariant $2$,
and a $2$-isogeny excludes them.

The groups $\torz{51}$, $\torz{63}$ and $\torz{65}$ (\Cref{sec:formal}) required a new
idea. For these $n$, every elliptic curve with a point of order $n$ over a sextic field has bad
reduction at all primes above $2$, but \cite[Proposition 5.2]{DN} does not apply: for $n=51$
our lower bound for the gonality of $X_1(51)$ is too weak, and for $n=63$ and $65$ there are
newforms of level dividing $n$ with positive analytic rank, so we cannot use that $J_1(n)(\Q)$
is finite. We show that the cuspidal reduction cannot occur by a formal
immersion argument at $2$. The usual difficulty at $p=2$, namely that reduction is not
injective on the $2$-torsion of the Jacobian, is resolved by \Cref{lem:dyadic}, which shows
that $T_2-\dm 2$ kills the $2$-torsion in the kernel of reduction, and the formal immersion
condition then becomes a condition on the Cartier operator on the holomorphic differentials of
a quotient of $X_1(n)$ over $\F_2$, which can be checked by linear algebra over $\F_2$. For
$n=65$ (and $n=57$) the Jacobian has positive rank, and we first show that the elliptic curve
cannot have good reduction at a prime above $2$, using the known quadratic and cubic points on
$X_0(57)$ and $X_0(65)$ (\Cref{sec:57_65}). The group $\torz{57}$ is then eliminated by a
computation of the rational $2$-torsion of $J_1(57)$.

\Cref{sec:sporadic} contains the four sporadic groups and the proof of \Cref{thm:main}, and
\Cref{sec:computations} describes the computations.

\subsection*{Data availability and reproducibility}
All the code and all the logs of the computations described in this paper are available at
\begin{center}\url{https://github.com/F-Najman/sextic_torsion}\end{center}
The \texttt{README} there lists, statement by statement, which script establishes each claim
and which log it produced. In the digital form of this article, the words reporting a
computation, such as ``we \gitlink{verify_37.m}{compute}'', are hyperlinks to the script that
carries it out.

\subsection*{AI disclosure and comments}
AI models were extensively used in the work leading to this paper. After completing \cite{DN} (with Derickx) and \cite{N5} and giving talks about this work, I was asked several times whether we could do such a classification for number fields of degree $d>5$. My answer was that I thought it was not possible as it would require explicit computations that did not seem feasible to me. To check whether this was correct, in mid-October 2026 I asked AI whether this was possible, and whether it could do it. To my surprise, it was able to complete the classifications of possible torsion groups of elliptic curves over number fields, not only in degree 6, but later also in degrees 7, 8, 9 and finally 10. After some initial hesitation about releasing primarily AI-generated mathematics, I decided that the results were simply too interesting and the mathematics too beautiful to just stay on my computer. 

OpenAI's GPT-6.1 Sol produced almost all of the mathematical content of the first version of the proof and carried out most of the computations. Anthropic's Claude Opus 5.5 reviewed that version, independently recomputed the relevant computations and the explicit examples, wrote the draft of the present manuscript and organized the code repository. I carefully reviewed the paper line by line and rewrote most parts of it. I understand it completely, and take full responsibility for the content of the paper.

\subsection*{Future work}

As noted in the AI disclosure above, I have obtained AI-generated proofs that appear to prove a complete classification of the possible torsion groups of elliptic curves over number fields of degrees 7, 8, 9 and 10. Once I have fully understood and verified the arguments in degrees 7--9, I plan to combine the results for degrees 6--9 into a single paper, which I currently expect to be around 120 pages long. Degree 10 will stay in a separate paper for mathematical reasons and because of its length.


\section{Notation and auxiliary results}\label{sec:prelim}

Throughout, $m\mid n$ are positive integers with $mn>4$. We
denote by $X:=X_1(m,n)$ the modular curve parametrizing triples $(E,P,Q)$ consisting of a
generalized elliptic curve $E$ and points $P,Q$ of orders $m$ and $n$ generating a subgroup of
$E[n]$ isomorphic to $\tg{m}{n}$, such that the Weil pairing of $P$ and $(n/m)Q$ is $\zeta_m$. It is a smooth
projective geometrically irreducible curve over $\Z[\zeta_m,\tfrac1n]$, and it is a fine moduli
space when $n>4$ (see \cite[\S 2]{DN} and \cite[\S 2]{DS}). Every curve occurring below has
$n\geq 5$, so we use the fine moduli interpretation throughout. We write $Y_1(m,n)$ for the open
subscheme obtained by removing the cusps, and $J:=J_1(m,n)$ for the Jacobian of $X$. More
generally, we write $J(Y)$ for the Jacobian of a curve $Y$, and $J$ for $J(X)$ when the curve $X$
is clear from the context. We abbreviate $X_1(n):=X_1(1,n)$ and $J_1(n):=J_1(1,n)$. For $m\leq 2$ the curve $X$ is defined
over $\Q$. For $d\geq1$ we write $X^{(d)}$ for the $d$-th symmetric power of $X$, so that
$X^{(d)}(\Q(\zeta_m))$ is the set of effective $\Q(\zeta_m)$-rational divisors of degree $d$ on $X$. For a curve
$X$ over a field $F$ with a set of cusps, such as $X_1(m,n)$ over $\Q(\zeta_m)$ or a quotient of it, we
write $\ClCusp_F X$ for the subgroup of $J(F)$ generated by the classes of the $F$-rational
divisors of degree $0$ supported on the cusps. As sums and negatives of such divisors are again
such divisors, every element of $\ClCusp_F X$ is the class of a single (not
necessarily effective) $F$-rational cuspidal divisor of degree $0$. For a divisor $D$ on a curve $X$ over a field $k$
we write $L(D):=\{f\in k(X)^\times:\ddiv(f)+D\geq0\}\cup\{0\}$ for its Riemann--Roch space.
Similarly, we write $\Omega(D)$ for the space consisting of $0$ and the non-zero differentials
$\omega$ on $X$ with $\ddiv(\omega)\geq D$. Then $\Omega(0)=H^0(X,\Omega^1)$ is the space of
holomorphic differentials on $X$, and if $D=\sum_P m_PP$ is effective, then $\Omega(D)$ is its
subspace of the differentials that vanish to order at least $m_P$ at every $P$. The
Riemann--Roch theorem states that $\dim L(D)=\deg D+1-g+\dim\Omega(D)$, where $g$ is the genus
of $X$.

Let $K$ be a number field of degree $6$ and $E/K$ an elliptic curve with
$\tg{m}{n}\subseteq E(K)$. By the Weil pairing $\Q(\zeta_m)\subseteq K$, so $\varphi(m)\mid 6$, and we
put $d_m:=6/\varphi(m)=[K:\Q(\zeta_m)]$. Choosing generators $P,Q$ with the right Weil pairing we obtain
a point $x\in Y_1(m,n)(K)$. We write $x^{(d_m)}$ for the sum of the
$\Q(\zeta_m)$-conjugates of $x$, each counted $[K:\Q(\zeta_m)(x)]$ times. It is an effective $\Q(\zeta_m)$-rational
divisor of degree $d_m$ on $X$ whose support contains no cusp, and we call it the \emph{trace
divisor} of $x$. If $m\leq2$, then $\Q(\zeta_m)=\Q$ and $d_m=6$, and we also write $x^{(6)}$.
All the arguments below are applied to trace divisors, so they also cover the
case in which the residue field $\Q(\zeta_m)(x)$ is smaller than $K$. We say that two pairs
$(K,E)$ and $(K',E')$, consisting of a number field and an elliptic curve over it, are
isomorphic if there is a field isomorphism $\sigma:K\to K'$ such that $E^\sigma\simeq E'$ over
$K'$.

The cusps of $X_1(m,n)$ correspond to N\'eron $k$-gons with $m\mid k\mid n$ together with a level
structure (see \cite[\S 2]{DN}). All cusps are defined over $\Q(\zeta_n)$, and
$\Gal(\Q(\zeta_n)/\Q(\zeta_m))$ acts on them through its action on $\mathbb G_m$. The cuspidal
subscheme is finite \'etale over $\Z[\zeta_m,\tfrac1n]$. We denote by $c_0$ the cusp given by
the N\'eron $n$-gon together with the level structure $(a,b)\mapsto (\zeta_m^a,b)$ (see
\cite[Lemma 2.8]{DN}), which is $\Q(\zeta_m)$-rational. From the description of the cusps
of $X_1(n)$ in \cite[\S 2]{DN} it follows that a cusp of $X_1(n)$ lying on a N\'eron $d$-gon is
fixed by the diamond operator $\dm a$ (defined below) only if $a\equiv\varepsilon\pmod d$ and
$a\equiv\varepsilon\pmod{n/d}$ for a common sign $\varepsilon=\pm1$. In particular $a$ is then
congruent to $\varepsilon$ modulo every prime dividing $n$.

For a congruence subgroup $\Gamma\subseteq\SL_2(\Z)$ we write $\overline\Gamma$ for its image
in $\PSL_2(\Z)$. The curve $X_1(m,n)_\C$ is the modular curve of the group
$\Gamma_1(n)\cap\Gamma(m)$, and for $n>2$ we have
$[\PSL_2(\Z):\overline{\Gamma_1(n)\cap\Gamma(m)}]=\tfrac m2\,n^2\prod_{\ell\mid n}(1-\ell^{-2})$,
the product being over the primes dividing $n$.

For $a\in(\torz n)^\times$ we write $\dm a$ for the diamond automorphism
$(E,P,Q)\mapsto (E,a^{-1}P,aQ)$ of $X_1(m,n)$, which preserves the Weil pairing. For $m\leq 2$
it is the usual diamond operator $(E,P,Q)\mapsto (E,aP,aQ)$. It depends only on $a$ modulo
$\pm1$, and we call $\Delta(n):=(\torz n)^\times/\{\pm1\}$ the group of diamond operators. We
write $\Z[\Delta(n)]$ for its group ring. It acts on divisors, on divisor classes and on $J$.
For $a_1,\dots,a_r\in(\torz n)^\times$ we write $X/\dm{a_1,\dots,a_r}$ for the quotient of $X$ by
the subgroup of $\Delta(n)$ generated by $\dm{a_1},\dots,\dm{a_r}$.
For a prime $q\nmid n$ with $q\equiv 1\pmod m$ we write $T_q$ for the Hecke correspondence
$T_q(E,P,Q)=\sum_{G}(E/G,P\bmod G,Q\bmod G)$, the sum being over the $q+1$ subgroups
$G\subseteq E[q]$ of order $q$, and we set
\[
A_q:=T_q-q\dm q-1 .
\]
The congruence $q\equiv1\pmod m$ ensures that $T_q$ preserves the Weil pairing, and that $q$
splits completely in $\Q(\zeta_m)$. The operator $A_q$ commutes with $\Z[\Delta(n)]$.

\begin{proposition}\label{prop:tools}
Let $q\nmid 2n$ be a prime with $q\equiv1\pmod m$.
\begin{enumerate}
\item \textup{(Eichler--Shimura)} $T_{q,\F_q}=\Frob_q+\dm q_{\F_q,*}\Ver_q$ on $J_{\F_q}$, where
      $\F_q$ is the residue field of a prime of $\Q(\zeta_m)$ above $q$.
\item If $L\supseteq \Q(\zeta_m)$ is a number field in which $q$ splits completely then
      $A_q\bigl(J(L)_\tors\bigr)=0$.
\item $A_q(c_0)=0$ as a divisor.
\item If $\rk J(\Q(\zeta_m))=0$ then $[A_q(D)]=0$ for every $D\in X^{(d)}(\Q(\zeta_m))$ and every $d$.
\end{enumerate}
\end{proposition}

\begin{proof}
For $m\leq 2$ these are, respectively, \cite[(2.2)]{DN}, \cite[Proposition 2.6]{DN},
\cite[Lemma 2.8]{DN} and \cite[Lemma 5.4]{DN}. For $m\geq 3$ the curve $X_1(m,n)$ is, over
$\Q(\zeta_m)$, isomorphic to the quotient of $X_1(mn)$ by the diamond operators congruent to $1$ modulo
$n$ \cite[\S 2]{DS}, and this isomorphism is compatible with $T_q$ and $\dm q$ for
$q\equiv1\pmod m$. Hence (1) follows from the Eichler--Shimura relation on $X_1(mn)$
\cite[(10.2.3)]{DiamondIm}, and the proofs of \cite[Proposition 2.6, Lemma 2.8 and Lemma 5.4]{DN}
carry over verbatim.
\end{proof}

We will also use the following results.

\begin{itemize}
\item (\emph{Manin--Drinfeld \cite{manin,drinfeld}}) the class of a degree $0$ cuspidal divisor
      on $X$ is torsion in $J(\Qbar)$;
\item (\emph{Katz \cite[Appendix]{katz}}) if $p>2$ is a prime of good reduction of $J$ and
      $\wp$ is a prime of a number field $L$ above $p$ which is unramified in $L$, then
      reduction modulo $\wp$ is injective on $J(L)_\tors$;
\item (\emph{Parent \cite[Proposition 2.4]{parent}}) if $A$ is an abelian variety over a
      number field $L$ with good reduction at a prime $\lambda$ above $2$ which is unramified
      in $L$, then the kernel of reduction modulo $\lambda$ on $A(L)_\tors$ is trivial or of
      exponent $2$. Indeed, this kernel is contained in the group of points of the formal group
      of the N\'eron model of $A$ at $\lambda$ with coordinates in the maximal ideal, and
      multiplication by $2$ maps it into the subgroup of points with coordinates in the square
      of the maximal ideal, which is torsion-free as $\lambda$ is unramified over $2$. We
      apply this to $J$ for $2\nmid n$ and to the Jacobians of quotients of $X_1(n)$ with good
      reduction at $2$;
\item (\emph{Abramovich \cite[Theorem 0.1]{abramovich}, with the bound
      $\lambda_1\geq 975/4096$ of Kim--Sarnak \cite[Appendix 2]{kim}}) for a congruence subgroup
      $\Gamma\subseteq\PSL_2(\Z)$,
      \[
        \gon_\C X_\Gamma\geq\tfrac{325}{2^{15}}\,[\PSL_2(\Z):\Gamma];
      \]
\item (\emph{Kato \cite[Corollary 14.3]{kato}}) if $A$ is a quotient of $J_1(n)$ over $\Q$ with
      $L(A,1)\neq0$, then $A(\Q)$ is finite;
\item (\emph{Clark--Genao--Pollack--Saia \cite{CGPS}}) the exact values of the least degree
      $d_{\mathrm{CM}}X_1(n)$ of a CM point on $X_1(n)$ that we quote are taken from the list
      \texttt{dcm\_list\_all\_min\_orders\_X1\_10k.m} in the repository
      \url{https://github.com/fsaia/least-cm-degree} accompanying \cite{CGPS}.
\end{itemize}
We \gitlink{global_table.py}{compute} all the values of Abramovich's bound and all the orders of
diamond operators that we quote below.

Finally, we state the two standard reduction lemmas.

\begin{lemma}\label{lem:additive}
Let $K$ be a number field, let $E/K$ be an elliptic curve with $\tg{m}{n}\subseteq E(K)$, and
let $\wp$ be a prime of $K$ above a rational prime $p\nmid n$. Then $E$ does not have additive
reduction at $\wp$.
\end{lemma}

\begin{proof}
This is contained in the proof of \cite[Proposition 5.2]{DN}. If $E$ had additive reduction at
$\wp$, the special fibre of the N\'eron model would have $a\cdot\#\F_\wp$ points over $\F_\wp$,
where $a\leq4$ is the order of the component group and $\#\F_\wp$ is a power of $p$. Reduction is
injective on the prime-to-$p$ torsion of $E(K)$, so $mn\mid a\leq 4$, contradicting $mn>4$.
\end{proof}

\begin{lemma}[{\cite[Lemma 5.3]{DN}}]\label{lem:badred}
Let $p\nmid n$ be a prime and let $K\supseteq \Q(\zeta_m)$ be a number field. Suppose that for every
prime $\wp$ of $K$ above $p$ there is no elliptic curve $E'$ over the residue field $\F_\wp$
with $\tg{m}{n}\subseteq E'(\F_\wp)$. Then no elliptic curve $E/K$ with
$\tg{m}{n}\subseteq E(K)$ has good reduction at a prime of $K$ above $p$.
\end{lemma}

Whenever we invoke \Cref{lem:badred}, the hypothesis follows from the Hasse bound, from
Waterhouse's description \cite[Theorem 4.1]{waterhouse} of the possible traces of Frobenius,
or from an explicit \gitlink{torsion_over_Fq.m}{enumeration} of all elliptic curves over the
relevant finite fields. For example, over $\F_{64}$ the possible traces of Frobenius are the odd
integers $t$ with $|t|\leq 16$ together with $t\in\{0,\pm8,\pm16\}$, so no elliptic curve
over $\F_{64}$ has exactly $51$, $55$, $63$, $75$ or $77$ points. If $\#E(\F_{64})=81$, then
the Frobenius $\pi$ satisfies $(\pi+8)^2=0$, so $\pi=-8$ as $\End(E)$ has no non-zero nilpotents,
and $E(\F_{64})=\ker(1-\pi)=E[9]$. Similarly, if $\#E(\F_{64})=49$, then $\pi=8$ and
$E(\F_{64})=E[7]$.

We will often use the following simple lemma. Its first part is the argument at the end of the
proof of \cite[Proposition 7.2]{N5}.

\begin{lemma}\label{lem:orbits}
Let $H\subseteq\Delta(n)$ be a subgroup and let $D$ be a divisor on $X_1(n)$ whose support
consists of non-cuspidal non-CM points. Then no element $h\neq1$ of $H$ fixes a geometric point in
the support of $D$. If moreover $hD=D$ for all $h\in H$, then the geometric points in the support of
$D$ form a union of $H$-orbits of length $\#H$, and $D=\sum_O m_O\sum_{x\in O}x$, where $O$ runs
over these orbits. In particular $\#H\mid\deg D$.
\end{lemma}

\begin{proof}
If $\dm a$ fixes a geometric point $x=(E,P)$ in the support of $D$, then $(E,aP)\simeq(E,P)$, so
there is an automorphism $u$ of $E$ with $u(P)=aP$. As $E$ has no complex multiplication,
$\Aut(E)=\{\pm1\}$, so $aP=\pm P$, and as $P$ has order $n$, we get $a\equiv\pm1\pmod n$, that
is, $\dm a=1$ in $\Delta(n)$. If $hD=D$ for all $h\in H$, then $D$ has the same coefficient at
$x$ and at $hx$, so the geometric points in the support of $D$ form a union of $H$-orbits and the
coefficients of $D$ are constant on each of them. By the first part, each orbit has length $\#H$.
\end{proof}

\section{The groups that have to be considered}\label{sec:groups}

Call a finite abelian group $G$ a \emph{sextic torsion group} if $G\simeq E(K)_\tors$ for some
elliptic curve $E$ over some number field $K$ of degree $6$. Every such $G$ is isomorphic to
$\tg{m}{n}$ with $m\mid n$. Put
\[
\mathcal B:=\Phi^\infty(6)\cup\{\torz{25},\torz{37},\tg{4}{12},\tg{2}{24}\}.
\]
We will show in \Cref{sec:sporadic} that all the groups in $\mathcal B$ are sextic torsion
groups. Note that $\mathcal B$ is closed under passing to subgroups.

\begin{lemma}\label{lem:66}
Let $\mathcal{G}$ be the set of groups $G=\tg{m}{n}$, $m\mid n$, such that
\begin{enumerate}
\item $G\notin\mathcal B$, while $H\in\mathcal B$ for every proper subgroup $H\subsetneq G$;
\item $\varphi(m)\mid 6$;
\item every prime divisor of $n$ lies in $S(6)=\{2,3,5,7,11,13,17,19,37\}$.
\end{enumerate}
Then $\mathcal{G}$ consists of the following $66$ groups:
\[
\begin{array}{l}
\torz{n}\ \text{ for }\ n\in\{32,33,34,35,36,38,39,40,42,44,45,48,49,50,\\
\qquad\quad 51,52,54,55,56,57,60,63,65,75,77,81,85,91,95,\\
\qquad\quad 119,121,125,133,143,169,187,209,221,247,289,323,361\},\\[2pt]
\torz{37n}\ \text{ for }\ n\in\{2,3,5,7,11,13,17,19,37\},\\[2pt]
\tg{2}{2n}\ \text{ for }\ n\in\{11,13,14,15\},\\[2pt]
\tg{3}{3n}\ \text{ for }\ n\in\{5,6,7,8,9\},\\[2pt]
\tg{4}{4n}\ \text{ for }\ n\in\{4,5,6\},\\[2pt]
\tg{6}{12},\quad \tg{7}{7},\quad \tg{9}{9}.
\end{array}
\]
Moreover every sextic torsion group not in $\mathcal B$ contains a member of $\mathcal G$.
\end{lemma}

\begin{proof}
Conditions (2) and (3) are necessary for $G$ to be a subgroup of $E(K)$ for a sextic field $K$.
The Weil pairing forces $\Q(\zeta_m)\subseteq K$, so $\varphi(m)\mid 6$, and if $p\mid n$ then
there is a point of order $p$ over a sextic field, so $p\in S(6)$ by \cite[Theorem 1.2]{DKSS}.
Since $\mathcal B$ is closed under passing to subgroups, a sextic torsion group not in
$\mathcal B$ contains a minimal group not in $\mathcal B$, which satisfies (1), (2) and (3).
The set $\mathcal G$ is finite: if $G=\torz n$ is not of prime order, then $\Z/(n/p)\Z\in\mathcal B$
for every prime $p\mid n$, so $n/p\leq 37$ and $n\leq 37^2$, and if $m>1$ then
$\torz n\in\mathcal B$, so $n\leq 37$. We \gitlink{group_list.py}{compute} the list by a finite
enumeration.
\end{proof}

Let $G=\tg{m}{n}\in\mathcal G$. An elliptic curve $E$ over a sextic field $K$ with
$G\hookrightarrow E(K)$ gives a point $x\in Y_1(m,n)(K)$ and a trace divisor $x^{(d_m)}$ as in
\Cref{sec:prelim}. Conversely, a point of $Y_1(m,n)$ defined over a sextic field gives such a
pair $(K,E)$, as $X_1(m,n)$ is a fine moduli space. Hence $G$ is a subgroup of the torsion of an
elliptic curve over a sextic field if and only if $X_1(m,n)$ has a non-cuspidal point
whose residue field is contained in a field of degree $d_m$ over $\Q(\zeta_m)$. To rule out $G$ we
show that no effective $\Q(\zeta_m)$-rational divisor of degree $d_m$ on $X_1(m,n)$ with non-cuspidal
support can be a trace divisor.

Several members of $\mathcal G$ are reduced to others by a $2$-isogeny.

\begin{lemma}\label{lem:4n}
Let $k\geq 2$. There is a $\Q$-rational morphism
$X_1(4k)\to X_1(2,2k)$, $(E,P)\mapsto \bigl(E/\dm{2kP},\,Q,\,P\bmod \dm{2kP}\bigr)$, where $Q$
generates $E[2]/\dm{2kP}$. Hence, if $\tg{2}{2k}$ is not a subgroup of the torsion of an
elliptic curve over a sextic field, then neither is $\torz{4k}$.
\end{lemma}

\begin{proof}
The morphism is the one described in \cite[\S 2]{DN} (see also \cite[\S 2]{DS}). It maps
non-cuspidal points to non-cuspidal points, and the image of a point defined over $K$ is
defined over $K$.
\end{proof}

We will apply \Cref{lem:4n} with $k=11,13,14,15$ to rule out $\torz{44}$, $\torz{52}$,
$\torz{56}$ and $\torz{60}$. For $k=12$ the hypothesis fails, as $\tg{2}{24}$ is a sextic
torsion group, and $\torz{48}$ requires a separate argument, given in \Cref{prop:48}.

The rest of the paper treats the $66$ groups of \Cref{lem:66} one family at a time.
\Cref{table:summary} shows where each of them is dealt with.

\begin{table}[H]
\centering
\footnotesize
\begin{tabular}{@{}l@{\quad}l@{}}
\toprule
$G\in\mathcal G$ & treated in\\
\midrule
$\torz{37n}$, $n=2,3,5,7,11,13,17,19,37$ & \Cref{prop:37}\\
$\torz n$, $n=119,121,125,143,169,187,$ & \Cref{prop:global} (global method)\\
\quad $209,221,289,323,361$ & \\
$\torz n$, $n=77,85,91,95,133,247$ & \Cref{prop:large} (diamond operators)\\
$\torz n$, $n=55,75,81$, and $\tg{3}{27}$, & \Cref{prop:rank0_cases}\\
\quad $\tg{4}{24}$, $\tg{7}{7}$, $\tg{9}{9}$ & \\
$\torz{45}$, $\tg{3}{15}$, $\tg{3}{21}$ & \Cref{prop:doubled_cases}\\
$\torz{54}$ & \Cref{prop:54}\\
$\torz n$, $n=33,34,35,38,39,40,42,49,50$ & \Cref{prop:sieve,prop:49}\\
$\tg{2}{22}$, $\tg{2}{26}$, $\tg{2}{28}$ & \Cref{prop:sieve,prop:2_26}\\
$\tg{3}{24}$, $\tg{4}{16}$, $\tg{4}{20}$, $\tg{6}{12}$ & \Cref{prop:sieve_m}\\
$\torz{32}$, $\torz{36}$, $\tg{3}{18}$, $\tg{2}{30}$ & \Cref{prop:32,prop:36,prop:3_18,prop:2_30}\\
$\torz{51}$, $\torz{63}$ & \Cref{prop:51_63} (formal immersion at $2$)\\
$\torz{57}$, $\torz{65}$ & \Cref{prop:57,prop:65}\\
$\torz n$, $n=44,52,56,60$ & \Cref{lem:4n}\\
$\torz{48}$ & \Cref{prop:48}\\
\bottomrule
\end{tabular}
\medskip
\caption{Treatment of the $66$ groups of \Cref{lem:66}.}
\label{table:summary}
\end{table}

\section{The global method and diamond operators}\label{sec:global}

In this section we eliminate the seventeen groups $\torz n$ with $n\geq 77$, $37\nmid n$ and
$n\neq81$.
Throughout the section we use the following consequence of Kato's theorem, as in
\cite[\S 7]{N5}. If $\psi\in\Z[\Delta(n)]$ annihilates every isogeny factor $A_{[f]}$ of
$J_1(n)$ with $L(f,1)=0$, then $\psi\bigl(J_1(n)(\Q)\bigr)$ is finite. This holds because
$J_1(n)$ is isogenous to $\bigoplus_{n'\mid n}\bigoplus_{[f]}A_{[f]}^{\sigma_0(n/n')}$, the inner sum being
over the Galois orbits of newforms of weight $2$ and level $n'$, and this decomposition is
equivariant for the diamond operators, $\dm a$ acting on $A_{[f]}$ through $\chi_f(a)$, where
$\chi_f$ is the nebentypus of $f$. For an $f$ with $L(f,1)\neq0$ the group $A_{[f]}(\Q)$ is
finite by Kato's theorem. In particular, an element $\psi=\sum_a c_a\dm a$ has this property if
$\sum_a c_a\chi_f(a)=0$ for every newform $f$ of level dividing $n$ with $L(f,1)=0$. All the
analytic ranks that we use are taken from the LMFDB \cite{lmfdb}, where they are recorded as
proved, and we \gitlink{lmfdb_ranks.py}{check} them and the character values below.

\subsection{The global method}
The main tool is \cite[Proposition 4.3]{DN}, due to Derickx. We use it in the form of
\cite[Proposition 4.1]{N5}: let $n$ be an integer with $3\nmid n$, let $a\in(\torz n)^\times$
be such that $(\dm a-1)J_1(n)(\Q)$ is finite, and put $k_3:=7$ if $a\in\{3,3^{-1}\}$ and
$k_3:=8$ otherwise. If
\begin{equation}\label{eq:global}
6k_3<\frac{325}{2^{15}}[\PSL_2(\Z):\overline{\Gamma_1(n)}],\qquad \ord\dm a\nmid 6,
\qquad d_{\mathrm{CM}}X_1(n)>6,
\end{equation}
then $X_1(n)$ has no non-cuspidal point of degree dividing $6$. The proof of
\cite[Proposition 4.1]{N5} applies without change with $d=6$ in place of $d=5$. The point is that
\cite[Proposition 4.3]{DN} shows that an effective $\Q$-rational divisor of degree $6$ on
$X_1(n)$ without cusps in its support, such as a trace divisor, is a sum of orbits of $\dm a$, and
by \Cref{lem:orbits} each orbit has length $\ord\dm a$ when the support contains no CM point.

\begin{proposition}\label{prop:global}
The group $\torz n$ is not a subgroup of the torsion of an elliptic curve over a sextic field
for
\[
n\in\{119,\,121,\,125,\,143,\,169,\,187,\,209,\,221,\,289,\,323,\,361\}.
\]
\end{proposition}

\begin{proof}
We use the elements $a$ from \cite[Table 2]{N5}, which are listed in \Cref{table:global}.
For these $a$ the group $(\dm a-1)J_1(n)(\Q)$ is finite by the beginning of this section, as for
$n\in\{119,121,125,169,221,289,361\}$ every newform of weight $2$ and level dividing $n$
with positive analytic rank has trivial character, and for the remaining $n$ the newforms of
positive analytic rank and non-trivial character have character orbits
\lmfdcharorbit{143}{h}, \lmfdcharorbit{187}{g}, \lmfdcharorbit{209}{f} and
\lmfdcharorbit{323}{k}, and $a$ lies in the kernel of each of them. The remaining conditions
of \eqref{eq:global} are \gitlink{global_table.py}{checked} in \Cref{table:global}.
\end{proof}

\begin{table}[H]
\centering
\begin{tabular}{rcccccc}
\toprule
$n$ & $a \bmod n$ & $\ord\dm a$ & $i(n)$ & $k_3$ & $\tfrac{325}{2^{15}}i(n)$ & $d_{\mathrm{CM}}X_1(n)$\\
\midrule
$119$ & $3$   & $48$  & $6912$  & $7$ & $68.55$  & $96$\\
$121$ & $2$   & $55$  & $7260$  & $8$ & $72.00$  & $110$\\
$125$ & $3$   & $50$  & $7500$  & $7$ & $74.38$  & $50$\\
$143$ & $67$  & $12$  & $10080$ & $8$ & $99.97$  & $120$\\
$169$ & $3$   & $39$  & $14196$ & $7$ & $140.79$ & $52$\\
$187$ & $122$ & $16$  & $17280$ & $8$ & $171.38$ & $160$\\
$209$ & $34$  & $18$  & $21600$ & $8$ & $214.23$ & $180$\\
$221$ & $3$   & $48$  & $24192$ & $7$ & $239.94$ & $96$\\
$289$ & $3$   & $136$ & $41616$ & $7$ & $412.75$ & $136$\\
$323$ & $101$ & $18$  & $51840$ & $8$ & $514.16$ & $288$\\
$361$ & $3$   & $171$ & $64980$ & $7$ & $644.48$ & $114$\\
\bottomrule
\end{tabular}
\medskip
\caption{Data for \Cref{prop:global}. Here $i(n)=[\PSL_2(\Z):\overline{\Gamma_1(n)}]$, and the
sixth column is Abramovich's lower bound for $\gon_\C X_1(n)$, rounded down to two decimal places, and it is larger than $6k_3\leq48$ in
every row. The orders $\ord\dm a$ are computed in $\Delta(n)$, and none of them divides $6$.}
\label{table:global}
\end{table}

\subsection{Several diamond operators}
For the levels $n=77,85,91,95,133,247$ the global method does not apply in degree $6$, either
because $\ord\dm a\mid 6$ for every $a$ with $(\dm a-1)J_1(n)(\Q)$ finite, or because the
gonality bound is too weak. In \cite{N5} the levels $91$, $95$, $133$ and $247$ were handled by
the global method, with $a=3,11,64,144$, but the orders of these diamond operators, $6$, $3$, $3$
and $3$, divide $6$. For these levels we work at the prime $2$, and we use more general
elements of the group ring of the diamond operators, such as $1+\dm{31}$ and $(\dm6-1)(\dm3-1)$,
together with an argument with $X_0(91)$. For $\psi=\sum_a c_a\dm a\in\Z[\Delta(n)]$
we write $\|\psi\|:=\sum_a|c_a|$.

\begin{proposition}\label{prop:several}
Let $n$ be odd, let $\psi\in\Z[\Delta(n)]$ be such that $\psi\bigl(J_1(n)(\Q)\bigr)$ is finite,
and assume that $\gon_\C X_1(n)>12\|\psi\|$. Let $K$ be a sextic field, let $E/K$ be an
elliptic curve with a point of order $n$, and suppose that $E$ has bad reduction at every prime
of $K$ above $2$. Then the trace divisor $D=x^{(6)}$ satisfies $\psi(D)=0$.
\end{proposition}

\begin{proof}
By \Cref{lem:additive} the reduction of $E$ at every prime of $K$ above $2$ is multiplicative,
so, fixing a prime $\wp$ of $L:=\Q(\zeta_n)$ above $2$, the reduction of $D$ modulo $\wp$ is an
effective divisor of degree $6$ supported on the cusps. As the cuspidal subscheme is finite
\'etale over $\Z[\tfrac1n]$ and all cusps are defined over $L$, this reduction lifts to an
effective cuspidal divisor $C$ of degree $6$ on $X_1(n)_L$. Write
\[
\psi[D-C]=\psi[D-6c_0]+\psi[6c_0-C].
\]
The first summand lies in $\psi\bigl(J_1(n)(\Q)\bigr)$, which is finite, and the second is
torsion by Manin--Drinfeld. It follows that $\psi[D-C]\in J_1(n)(L)_\tors$, and it reduces to $0$
modulo $\wp$. As $2$ is unramified in $L$, Parent's result gives $2\psi[D-C]=0$, so
$2\psi(D-C)=\ddiv(f)$ for some non-zero rational function $f$ on $X_1(n)_L$. Writing $\psi=\psi^+-\psi^-$ with $\psi^\pm$
having non-negative coefficients, the polar divisor of $f$ is contained in
$2\psi^+(C)+2\psi^-(D)$, which has degree at most $12\|\psi\|$. By the assumption on the
gonality $f$ is constant, so $\psi(D)=\psi(C)$. The diamond operators preserve the set of cusps
and its complement, and $D$ has no cusps in its support, so $\psi(D)=0$.
\end{proof}

We will apply \Cref{prop:several} with $\psi=\dm a-1$, in which case its conclusion is that
$D$ is invariant under $\dm a$, and with a few other elements. Note that the case $\psi=\dm a-1$
is the degree $6$ version of \cite[Proposition 7.2]{N5}. If the support of $D$ contains no CM
point and $D$ is invariant under a subgroup $H\subseteq\Delta(n)$, then by \Cref{lem:orbits} it
is a sum of $H$-orbits of length $\#H$.

\begin{table}[H]
\centering
\begin{tabular}{rlrl}
\toprule
$n$ & characters of the newforms of positive analytic rank & $d_{\mathrm{CM}}X_1(n)$ & $\tfrac{325}{2^{15}}i(n)$\\
\midrule
$77$  & trivial & $60$ & $28.56$\\
$85$  & trivial, \lmfdcharorbit{85}{j} (order $4$) & $32$ & $34.27$\\
$91$  & trivial, \lmfdcharorbit{91}{g} (order $3$), \lmfdcharorbit{91}{u} (order $6$) & $24$ & $39.99$\\
$95$  & \lmfdcharorbit{95}{l} (order $12$) & $72$ & $42.84$\\
$133$ & trivial, \lmfdcharorbit{133}{f} (order $3$), \lmfdcharorbit{133}{s} (order $6$) & $36$ & $85.69$\\
$247$ & trivial, \lmfdcharorbit{247}{be} (order $12$) & $72$ & $299.92$\\
\bottomrule
\end{tabular}
\medskip
\caption{The newforms of weight $2$ and level dividing $n$ with positive analytic rank, listed
by their character orbits, and the data used in \Cref{prop:large}, with $i(n)$ as in
\Cref{table:global}. Their Conrey labels are $85.4$, $91.9$, $91.30$, $95.8$, $133.39$,
$133.31$ and $247.8$.}
\label{table:large}
\end{table}

\begin{proposition}\label{prop:large}
The group $\torz n$ is not a subgroup of the torsion of an elliptic curve over a sextic field
for $n\in\{77,85,91,95,133,247\}$.
\end{proposition}

\begin{proof}
Suppose $E$ is an elliptic curve over a sextic field $K$ with a point of order $n$. We first
show that $E$ has bad reduction at all primes of $K$ above $2$. By \Cref{lem:badred} it suffices
to show that no elliptic curve over $\F_{2^f}$, $f\leq 6$, has a point of order $n$. For $f\leq5$
this follows from the Hasse bound, since $(2^{5/2}+1)^2<45$, and for $f=6$ the Hasse bound gives
$\#E'(\F_{64})\leq 81$. This settles $n\neq 77$, and for $n=77$ we note that an elliptic curve
over $\F_{64}$ with exactly $77$ points would have trace of Frobenius $-12$, which is impossible
by \cite[Theorem 4.1]{waterhouse}, as $-12$ is even and not one of $0$, $\pm8$, $\pm16$.

Hence \Cref{prop:several} applies to $D=x^{(6)}$ for every $\psi$ as in its statement. We use
the data of \Cref{table:large}. In every case $d_{\mathrm{CM}}X_1(n)>6$, so the support of $D$
consists of non-CM points, and by \Cref{lem:orbits}, $D$ is a sum of orbits of length $\#H$ of
every subgroup $H$ of $\Delta(n)$ that fixes it.

For $n=77$ take $\psi=\dm3-1$. Only the trivial character occurs, so $\psi(J_1(77)(\Q))$ is
finite, and $12\|\psi\|=24<28.56$. Hence $\dm3D=D$, but $\ord\dm 3=30$ does not divide $6$,
a contradiction.

For $n=85$ take $\psi=\dm2-1$. The character of \lmfdcharorbit{85}{j} is trivial on $2$, and
$\ord\dm2=8\nmid 6$, so we conclude as for $n=77$.

For $n=91$ take $\psi=\dm3-1$. The characters of \lmfdcharorbit{91}{g} and
\lmfdcharorbit{91}{u} are trivial on $3$, so $\dm3D=D$. Here $\dm3$ has order $6$, so $D$ is a
single $\dm3$-orbit of length $6$, and in particular $x$ has degree $6$ and its $6$ conjugates
form one orbit under the diamond operators. It follows that the image of $x$ on
$X_0(91)=X_1(91)/\Delta(91)$ is a rational non-cuspidal point. But $X_0(91)$ has no rational
non-cuspidal points \cite{kenku91}, which is a contradiction.

For $n=95$ take $\psi=1+\dm{31}$. The only newform of positive analytic rank has character
$\chi$ with $\chi(31)=-1$, so $\psi$ annihilates it, and $12\|\psi\|=24<42.84$. Hence
$D+\dm{31}D=0$, which is impossible for a non-zero effective divisor.

For $n=133$, let $\chi_f$ and $\chi_s$ be characters in the orbits \lmfdcharorbit{133}{f}
and \lmfdcharorbit{133}{s}. We \gitlink{check_diamonds.m}{compute} that $\chi_f(6)=1$ and
$\chi_s(3)=1$, so
\[
\psi=(\dm6-1)(\dm3-1),\qquad \|\psi\|=4,\qquad 12\|\psi\|=48<85.69,
\]
annihilates every factor of positive analytic rank, and \Cref{prop:several} gives $\psi(D)=0$.
Hence the divisor $D'=\dm3D-D$ satisfies $\dm6D'=D'$. Its support consists of non-CM points, so
by \Cref{lem:orbits} its geometric points form a union of $\dm6$-orbits of length $\ord\dm6=18$.
But these are among the at most $12$ geometric points in the supports of $D$ and $\dm3D$, and
hence $D'=0$. So $\dm3D=D$, and
since $\dm3$ has order $9$, we get $9\mid\deg D=6$, which is a contradiction.

For $n=247$, the characters in \lmfdcharorbit{247}{be} are trivial on $12$ and on $56$, so
\Cref{prop:several} gives $\dm{12}D=\dm{56}D=D$. We \gitlink{check_diamonds.m}{compute} that
$\dm{12}$ and $\dm{56}$ generate a subgroup of order $9$ of $\Delta(247)$, so by \Cref{lem:orbits},
$D$ is a sum of orbits of length $9$ of this subgroup, and $9\mid\deg D=6$, which is a contradiction.
\end{proof}

\section{The rank zero criterion}\label{sec:rank0}

We recall \cite[Proposition 5.2]{DN} in the form in which we use it.

\begin{proposition}[{\cite[Proposition 5.2]{DN}}]\label{prop:main}
Let $m\mid n$ with $mn>4$, let $p\nmid mn$ be a prime and let $d\geq 1$. Assume
$\rk J_1(m,n)(\Q(\zeta_m))=0$ and that one of the following holds:
\begin{enumerate}
\item[\textup{(a)}] $p>2$ and $\gon_{\Q(\zeta_m)}X_1(m,n)>d$;
\item[\textup{(b)}] $p=2$ and $\gon_{\Q(\zeta_m)}X_1(m,n)>2d$.
\end{enumerate}
Then every elliptic curve $E$ over a degree $d$ extension $K$ of $\Q(\zeta_m)$ with
$\tg{m}{n}\subseteq E(K)$ has, for every prime $\mathfrak p$ of $\Q(\zeta_m)$ above $p$, good
reduction at at least one prime of $K$ above $\mathfrak p$.
\end{proposition}

In \cite[Proposition 5.2]{DN} the conclusion is stated for the primes of $K$ above $p$, but the
proof fixes a prime $\mathfrak p$ of $\Q(\zeta_m)$ above $p$ and uses only the reduction of $E$ at
the primes of $K$ above $\mathfrak p$, so it gives the statement above.

To apply this to $\tg{7}{7}$ and $\tg{9}{9}$ we need the following lemma, as these cases are
not covered by \cite[Theorem 4.1]{DS}.

\begin{lemma}\label{lem:fullrank}
For $m=7$ and $m=9$ we have $\rk J_1(m,m)(\Q(\zeta_m))=0$.
\end{lemma}

\begin{proof}
Put $F=\Q(\zeta_m)$. We use the argument of the proof of \cite[Theorem 4.1]{DS}, where $F$ is
quadratic: we show that the lemma follows from $L(g,1)\neq0$ for the newforms $g$ described below,
and then check this by a computation. Over $F$ the curve $X_1(m,m)$ is isomorphic to the quotient $X_H$ of
$X_1(m^2)$ by the diamond operators congruent to $1$ modulo $m$ \cite[\S 2]{DS}, so it suffices
to show that $B(F)$ is finite for the Jacobian $B$ of $X_H$, which is an abelian variety over
$\Q$. The space of cusp forms of weight $2$ on $X_H$ is
\[
V_m:=\bigoplus_{\varepsilon}S_2(\Gamma_0(m^2),\varepsilon),
\]
the sum being over the even Dirichlet characters $\varepsilon$ modulo $m$, so after extending
scalars to $\overline\Q_\ell$ the $\ell$-adic representation of $B$ is the sum of the
representations $\rho_f$ attached to the newforms $f$ occurring in $V_m$, with multiplicities.
We have $B(F)=R(\Q)$ for the restriction of scalars $R:=\operatorname{Res}_{F/\Q}B_F$, whose
$\ell$-adic representation is induced from the restriction of that of $B$ to
$\Gal(\overline\Q/F)$. Hence, over $\overline\Q_\ell$, it is the sum of the twists
$\rho_f\otimes\chi$, where $\chi$ runs over the characters of $\Gal(F/\Q)\simeq(\torz m)^\times$.
The twist of a newform of level dividing $m^2$ and character $\varepsilon$ modulo $m$ by a
character $\chi$ modulo $m$ is an eigenform of level dividing $m^2$ and character
$\varepsilon\chi^2$, which is again even and defined modulo $m$, so $\rho_f\otimes\chi\simeq\rho_g$
for a newform $g$ occurring in $V_m$. The multiset of these newforms $g$ is stable under the
action of $\Aut(\C)$ on the coefficients, so by Faltings' isogeny theorem \cite{faltings} the
abelian variety $R$ is isogenous over $\Q$ to a product of the abelian varieties $A_{[g]}$ attached
to the Galois orbits of the newforms $g$ occurring in $V_m$. By Kato's theorem it therefore suffices to show that
$L(g,1)\neq0$ for every newform $g$ occurring in $V_m$. Using modular symbols we \gitlink{full_level_rank.m}{compute}
the dimension of the winding submodule of each Galois orbit of newforms occurring in $V_m$, and
it is equal to the dimension of the orbit in every case. The newforms in question have level
$49$ for $m=7$, and levels $27$ and $81$ for $m=9$.
\end{proof}

\begin{proposition}\label{prop:rank0_cases}
The groups $\torz{55}$, $\torz{75}$, $\torz{81}$, $\tg{3}{27}$, $\tg{4}{24}$, $\tg{7}{7}$ and
$\tg{9}{9}$ are not subgroups of the torsion of an elliptic curve over a sextic field.
\end{proposition}

\begin{proof}
In each case suppose that $E$ is an elliptic curve over a sextic field $K$ containing the group
in question, and let $d=d_m=[K:\Q(\zeta_m)]$. We show that the hypotheses of \Cref{prop:main} hold for
a suitable prime $p$, and that $E$ cannot have good reduction at any prime of $K$ above $p$ by
\Cref{lem:badred}, which is a contradiction.

For $n=55,75,81$ we have $\rk J_1(n)(\Q)=0$ by \cite[Theorem 3.1(2)]{DEHMZ}, and Abramovich's
bound gives $\gon_\Q X_1(n)>14.28,23.80,28.92$, respectively, all larger than $12$, so
\Cref{prop:main}(b) applies with $p=2$ and $d=6$. No elliptic curve over $\F_{2^f}$ with
$f\leq 5$ has a point of order $n$ by the Hasse bound. Over $\F_{64}$ a point of order $55$ or
$75$ would force $\#E(\F_{64})=55$ or $75$, which is impossible, and a point of order $81$ would
force $E(\F_{64})=E[9]$, which has no point of order $81$ (see the discussion after \Cref{lem:badred}).

For $\tg{3}{27}$ we have $d=3$. By \cite[Theorem 4.1]{DS},
$\rk J_1(3,27)(\Q(\zeta_3))=0$, and the index $[\PSL_2(\Z):\overline{\Gamma_1(27)\cap\Gamma(3)}]=972$
gives $\gon_{\C}X_1(3,27)>9.64>6$, so \Cref{prop:main}(b) applies with $p=2$. The prime $2$ is
inert in $\Q(\zeta_3)$, so the residue fields of $K$ at the primes above $2$ are $\F_4$, $\F_{16}$ or
$\F_{64}$. Over the first two the Hasse bound excludes a subgroup of order $81$, and over
$\F_{64}$ a subgroup $\tg{3}{27}$ would force $E(\F_{64})=E[9]\simeq(\torz9)^2$ (see the discussion after \Cref{lem:badred}),
which has no point of order $27$.

For $\tg{4}{24}$ we have $d=3$. By \cite[Theorem 4.1]{DS},
$\rk J_1(4,24)(\Q(i))=0$, and the index $[\PSL_2(\Z):\overline{\Gamma_1(24)\cap\Gamma(4)}]=768$ gives
$\gon_\C X_1(4,24)>7.61>3$, so
\Cref{prop:main}(a) applies with $p=5$. The prime $5$ splits in $\Q(i)$, so the residue fields of
$K$ at the primes above $5$ are $\F_5$, $\F_{25}$ or $\F_{125}$. The Hasse intervals of these
fields contain no multiple of $96$, since $(\sqrt{25}+1)^2<96$ and
$96<(\sqrt{125}-1)^2<(\sqrt{125}+1)^2<192$.

For $\tg{7}{7}$ and $\tg{9}{9}$ we have $K=\Q(\zeta_m)$ and $d=1$. The genera of $X_1(7,7)$ and
$X_1(9,9)$ are $3$ and $10$, so $\gon_{\Q(\zeta_m)}X_1(m,m)>1$, and \Cref{lem:fullrank} gives
$\rk J_1(m,m)(\Q(\zeta_m))=0$. Take $p=29$ for $m=7$ and $p=19$ for $m=9$. Then $p\equiv1\pmod m$, so
the residue fields of $K$ at the primes above $p$ are $\F_p$, and since $(\sqrt p+1)^2<m^2$ no
elliptic curve over $\F_p$ has full $m$-torsion.
\end{proof}

\subsection{A criterion using twice a cuspidal divisor}
When $p=2$, \Cref{prop:main} requires $\gon_{\Q(\zeta_m)}X_1(m,n)>2d$, because the kernel of reduction
on the torsion of the Jacobian can contain elements of order $2$. For $X_1(45)$ and
$X_1(3,21)$ the gonality bounds that we have are too weak for this, and we replace the
gonality condition by a computation.

\begin{lemma}\label{lem:doubled}
Let $m\mid n$ with $n$ odd and $mn>4$, let $\mathfrak p$ be a prime of $\Q(\zeta_m)$ above $2$ with
residue field $k$, and let $d\geq1$. Put $X=X_1(m,n)$ and assume that $\rk J(\Q(\zeta_m))=0$ and that
\[
\dim L(2\overline C)=1
\]
for every effective cuspidal divisor $\overline C$ of degree $d$ on $X_k$. Then every elliptic
curve $E$ over a degree $d$ extension $K$ of $\Q(\zeta_m)$ with $\tg{m}{n}\subseteq E(K)$ has good
reduction at at least one prime of $K$ above $\mathfrak p$.
\end{lemma}

\begin{proof}
Suppose the opposite, and let $D$ be the trace divisor of the corresponding point. By
\Cref{lem:additive} the reduction of $E$ at every prime of $K$ above $\mathfrak p$ is
multiplicative, so $D$ reduces to an effective cuspidal divisor $\overline C$ of degree $d$ on
$X_k$. Let $M=\Q(\zeta_n)$, choose a prime $\mathfrak P$ of $M$ above $\mathfrak p$, and let
$F\subseteq M$ be the fixed field of its decomposition group in $\Gal(M/\Q(\zeta_m))$. The prime
$\mathfrak P\cap F$ of $F$ has residue field $k$, and $F$ is unramified at $2$ as $n$ is odd.
All cusps are defined over $M$ and the cuspidal subscheme is finite \'etale at $2$, so the
closed points in the support of $\overline C$ lift to orbits of $\Gal(M/F)$ of the same size,
and $\overline C$ lifts to an effective $F$-rational cuspidal divisor $C$ of degree $d$. The
class $[D-C]=[D-dc_0]+[dc_0-C]$ lies in $J(F)$ and is torsion, by the rank assumption and by
Manin--Drinfeld, and it reduces to $0$ modulo $\mathfrak P\cap F$. By Parent's result
$2[D-C]=0$, so $2D$ and $2C$ are linearly equivalent. By semicontinuity \cite[Theorem III.12.8]{hartshorne}, the dimension of the
Riemann--Roch space of $2C$ on $X_F$ satisfies
\[
\dim L(2C)\leq\dim L(2\overline C)=1,
\]
so $2C$ is the only effective divisor in its class and $2D=2C$, which is impossible as $D$ is
non-cuspidal and $C$ is cuspidal.
\end{proof}

We check the condition of \Cref{lem:doubled} using differentials. Let $g$ be the genus of
$X_k$. For a closed point $P$ of $X_k$ and $j\geq0$, the space $\Omega(jP)$ of holomorphic
differentials on $X_k$ that vanish to order at least $j$ at $P$ is the subspace of the
$g$-dimensional space $H^0(X_k,\Omega^1)$ defined by the vanishing of the first $j$ coefficients
of the expansion in a local parameter at $P$. These coefficients lie in the residue field of $P$,
so the codimension of $\Omega(jP)$ is at most $j\deg P$. For an effective divisor
$\overline C=\sum_P m_PP$, the space $\Omega(\overline C)$ is the intersection of the subspaces
$\Omega(m_PP)$, so it has codimension at most $\deg\overline C$. By the Riemann--Roch theorem,
\[
\dim L(2\overline C)=2d+1-g+\dim\Omega(2\overline C)=1+2d-\operatorname{codim}\Omega(2\overline C).
\]
As $\Omega(2\overline C)$ is cut out by $\sum_P2m_P\deg P=2d$ linear conditions over $k$ (the first
$2m_P$ coefficients at each $P$), we get $\dim L(2\overline C)=1$ if and only if these conditions
are linearly independent.

This is much more efficient than computing each Riemann--Roch space $L(2\overline C)$ directly.
The subspaces $\Omega(2mP)$ depend only on the closed cusp $P$ and on
$m\leq d/\deg P$, so each of them is computed only once, by expanding a basis of
$H^0(X_k,\Omega^1)$ in a local parameter at $P$. For $X_1(45)_{\F_2}$ these are $91$ subspaces,
while there are $19048$ divisors $\overline C$ to check. For each $\overline C$ it then remains to
intersect at most $6$ subspaces of $H^0(X_k,\Omega^1)$, which is fast linear algebra over $k$.

\begin{proposition}\label{prop:doubled_cases}
The groups $\torz{45}$, $\tg{3}{21}$ and $\tg{3}{15}$ are not subgroups of the torsion of an
elliptic curve over a sextic field.
\end{proposition}

\begin{proof}
Let $E$ be an elliptic curve over a sextic field $K$ containing one of these groups. In all
three cases $E$ has bad reduction at every prime of $K$ above $2$ by \Cref{lem:badred}. For
$\torz{45}$ the residue fields are $\F_{2^f}$ with $f\leq 6$, and no Hasse interval of these
fields contains a multiple of $45$. For $\tg{3}{21}$ and $\tg{3}{15}$ the residue fields are
$\F_4$, $\F_{16}$ and $\F_{64}$ as in the proof of \Cref{prop:rank0_cases}. An elliptic curve
over $\F_{64}$ with exactly $63$ points does not exist, and the Hasse intervals of these three
fields contain no other multiple of $63$ and no multiple of $45$. So in each case it suffices
to show that $E$ has good reduction at some prime above $2$.

For $\torz{45}$ we have $\rk J_1(45)(\Q)=0$ by \cite[Theorem 3.1(2)]{DEHMZ}, and we apply
\Cref{lem:doubled} with $m=1$, $d=6$ and $k=\F_2$. The curve $X_1(45)_{\F_2}$ has genus $41$,
and its closed cusps of degree at most $6$ have degrees $1^{12},2^4,4^5,6^2$ (where $a^r$ stands for $r$ cusps of degree $a$), so there are
$19048$ effective cuspidal divisors of degree $6$. For each of them we
\gitlink{cusp_double_rr_local_ref_lazy.m}{compute} that $\dim L(2\overline C)=1$.

For $\tg{3}{21}$ we have $\rk J_1(3,21)(\Q(\zeta_3))=0$ by \cite[Theorem 4.1]{DS}, and we apply
\Cref{lem:doubled} with $m=3$, $d=3$ and $k=\F_4$. The curve $X_1(3,21)_{\F_4}$ has genus $25$
and its cusps have degrees $1^{24},3^8$, so there are $\binom{26}{3}+8=2608$ effective
cuspidal divisors of degree $3$. For each of them we \gitlink{cusp_double_rr_3n.m}{compute}
that $\dim L(2\overline C)=1$.

For $\tg{3}{15}$ we argue similarly, with one extra step. Let $X=X_1(3,15)$ over $\Q(\zeta_3)$, a curve
of genus $9$ with $\rk J(\Q(\zeta_3))=0$ \cite[Theorem 4.1]{DS}, and let $\mathfrak p$ be the prime of
$\Q(\zeta_3)$ above $2$, with residue field $\F_4$. The curve $X_{\F_4}$ has $16$ cusps of degree $1$
and $8$ of degree $2$, so there are $\binom{18}{3}+16\cdot8=944$ effective cuspidal divisors of
degree $3$. We \gitlink{cusp_double_rr_3n_localref.m}{compute} that
$\dim L(2\overline C)=1$ for $912$ of them, so the proof of
\Cref{lem:doubled} excludes these as reductions of $D$. The remaining $32$ divisors are
supported on the cusps of degree $1$ and satisfy $\dim L(\overline C)=1$.
The $16$ cusps of $X_{\F_4}$ of degree $1$ are the reductions of the $16$ cusps of $X$ that are
$\Q(\zeta_3)$-rational. Indeed, for $a\in(\torz{15})^\times$ let $\sigma_a$ be the automorphism of
$\Q(\zeta_{15})$ with $\sigma_a(\zeta_{15})=\zeta_{15}^a$. Then $\Gal(\Q(\zeta_{15})/\Q(\zeta_3))$ is
cyclic of order $4$, generated by $\sigma_7$, and its Frobenius element at $\mathfrak p$ is
$\sigma_4=\sigma_7^2$. As the cuspidal subscheme is finite \'etale at $2$, a cusp of $X$ reduces to a
cusp of degree $1$ if and only if it is fixed by $\sigma_4$, and it is $\Q(\zeta_3)$-rational if and
only if it is fixed by $\sigma_7$. We \gitlink{cusp_orbits.py}{check} that each of $\sigma_7$ and
$\sigma_4$ fixes exactly $16$ of the $32$ cusps of $X$, so, as $\sigma_4=\sigma_7^2$, they fix the
same cusps. It follows that such a $\overline C$ lifts to a $\Q(\zeta_3)$-rational cuspidal divisor $C$,
and $[D-C]\in J(\Q(\zeta_3))$ reduces to $0$.

We claim that reduction modulo $\mathfrak p$ is injective on $J(\Q(\zeta_3))=J(\Q(\zeta_3))_\tors$. It is
injective on the odd part. Let $G\subseteq J(\Q(\zeta_3))$ be the subgroup generated by the
differences of the $16$ $\Q(\zeta_3)$-rational cusps. We \gitlink{verify_3_15_reduction_injective.m}{compute}
that the reduction of $G$ modulo $\mathfrak p$ has invariant factors $[3,3,24,96,96]$, whose
$2$-part has order $2^{13}$, and that $J(\F_7)$, for a prime of $\Q(\zeta_3)$ above $7$, has invariant
factors $[3,3,3,24,480,480]$, whose $2$-part also has order $2^{13}$. As reduction modulo a prime
above $7$ is injective on $J(\Q(\zeta_3))_\tors$, the $2$-part of $J(\Q(\zeta_3))$ has order at most $2^{13}$,
and as the reduction of $G$ has a $2$-part of order $2^{13}$, it follows that
$J(\Q(\zeta_3))[2^\infty]=G[2^\infty]$ and that reduction modulo $\mathfrak p$ is injective on it. This
proves the claim. It follows that $D\sim C$, and by semicontinuity \cite[Theorem III.12.8]{hartshorne}
$\dim L(C)\leq\dim L(\overline C)=1$, so $D=C$,
which is a contradiction.
\end{proof}

\section{The Hecke sieve in degree six}\label{sec:sieve}

Let $X=X_1(m,n)$ with $\rk J(\Q(\zeta_m))=0$, let $\mathfrak p$ be a prime of $\Q(\zeta_m)$ above a prime
$p\nmid n$, with residue field $k$, let $K\supseteq \Q(\zeta_m)$ be a field of degree $6$ and let
$x\in Y_1(m,n)(K)$, with trace divisor $D_x:=x^{(d_m)}$. Reducing modulo $\mathfrak p$ we may
write
\[
(D_x)_{k}=D+C,
\]
where $D$ is effective of degree $d'$ with non-cuspidal support and $C$ is effective cuspidal of
degree $d_m-d'$, and $d'\geq 1$ if and only if the elliptic curve corresponding to $x$ has good
reduction at some prime of $K$ above
$\mathfrak p$. As in \cite[\S 6]{N5} we call the multiset of degrees of the closed points
occurring in $D$, a partition of $d'$, the \emph{type} of the decomposition, and write it in
brackets. The Hecke sieve \cite[Proposition 5.5]{DN} says that if $[A_q(D)]\neq 0$ and
$[A_q(C)]=0$ in $J(k)$ for some prime $q\nmid 2pn$ with $q\equiv1\pmod m$, then this
decomposition cannot occur, since $[A_q(D_x)]=0$ by \Cref{prop:tools}(4) and its reduction is
$[A_q(D)]+[A_q(C)]$. Since $A_q$ commutes with the diamond operators, it suffices to test one $D$
in each orbit of the diamond action.

The condition $[A_q(C)]=0$ is supplied by \cite[Lemma 6.1]{N5}, whose proof applies to
$X_1(m,n)$ for every $m$. The cusps are defined over $M:=\Q(\zeta_n)$, and the Frobenius of
$\mathfrak p$ acts on them as the element $p^f$ of $\Gal(M/\Q(\zeta_m))\subseteq(\torz n)^\times$, where
$k=\F_{p^f}$. For a closed cusp $C_i$ of $X_k$ let $H_i\subseteq\Gal(M/\Q(\zeta_m))$ be the stabilizer of
the corresponding orbit of geometric cusps. Then $[A_q(C)]=0$ for every cuspidal divisor $C$
supported on closed cusps $C_i$ with $q\bmod n\in H_i$, and in particular whenever
$q\equiv p^{fj}\pmod n$ for some $j\geq0$. If a closed cusp $C_i$ is the reduction of a
$\Q(\zeta_m)$-rational cusp, then $H_i=\Gal(M/\Q(\zeta_m))$ and it imposes no condition on $q$.

\subsection{The curves over \texorpdfstring{$\Q$}{Q}}

\begin{proposition}\label{prop:sieve}
The groups $\torz n$ for $n\in\{33,34,35,38,39,40,42,50\}$ and the groups $\tg{2}{22}$ and
$\tg{2}{28}$ are not subgroups of the torsion of an elliptic curve over a sextic field. Hence,
by \Cref{lem:4n}, neither are $\torz{44}$ and $\torz{56}$.
\end{proposition}

\begin{proof}
Let $X=X_1(m,n)$ be one of the ten curves and let $p$ be the prime in \Cref{table:sieve}. In
every case $\rk J(\Q)=0$. For $m=1$ this is \cite[Theorem 3.1(2)]{DEHMZ}, and for $m=2$ it is
\cite[Theorem 3.1(3)]{DEHMZ} (or \cite[Theorem 4.1]{DS}).

We first exclude the totally cuspidal reduction using \Cref{prop:main} with $d=6$. This
requires $\gon_\Q X>6$ when $p>2$, and $\gon_\Q X>12$ when $p=2$, which happens only for $n=39$.
By \cite[Theorem 4 and Table 1]{DvH} we have $\gon_\Q X_1(n)=10,10,12,12,14,12$ for
$n=33,34,35,38,39,40$, and Derickx and van Hoeij proved $\gon_\Q X_1(42)>8$
\cite[proof of Proposition 6]{DvH}. Abramovich's bound gives $\gon_\Q X_1(50)>8.92$. For
$X_1(2,22)$ and $X_1(2,28)$, Derickx and Sutherland proved $\gon_{\F_3}X_1(2,2k)>6$
\cite[proof of Proposition 5.5]{DS}, so $\gon_\Q X_1(2,2k)>6$.

It remains to exclude every decomposition $(D_x)_{\F_p}=D+C$ with $\deg D=d'\geq 1$. The types
that can occur are limited by the degrees of the non-cuspidal closed points of $X_{\F_p}$. The
numbers of these points and of their diamond orbits are listed in \Cref{table:sieve}. For one
representative $D$ of each diamond orbit, and for every possible cuspidal complement $C$, we
compute, with the scripts linked in \Cref{table:sieve}, that $[A_q(D)]\neq0$ for one of the
listed primes $q$ for which $[A_q(C)]=0$ holds, and this leaves no survivors.

For example, $X_1(39)_{\F_2}$ has no non-cuspidal closed points of degree less than $6$, and
$12$ of degree $6$. As the totally cuspidal reduction is excluded, the reduction of the trace
divisor is one of these $12$ points, so it has type $[6]$ and $C=0$. The $12$ points form a
single orbit under $\Delta(39)$, so it suffices to take one of them as $D$, and we
\gitlink{hecke_sieve_deg5_v6.m}{compute} that $[A_q(D)]\neq0$ for $q=5$ or $q=7$.
\end{proof}

\begin{table}[H]
\centering
\small
\begin{tabular}{lrrlll}
\toprule
$X$ & $g$ & $p$ & non-cuspidal points of degree $1,\dots,6$ & orbits & $q$\\
\midrule
\gitlink{hecke_sieve_deg5_v6.m}{$X_1(32)$} & $17$ & $3$ & $0,0,8,8,16,140$          & $0,0,1,1,2,19$       & $5,7,11$\\
\gitlink{hecke_sieve_cached.m}{$X_1(33)$} & $21$ & $5$ & $0,10,60,230,602,2660$    & $0,1,6,24,61,269$    & $7$\\
\gitlink{hecke_sieve_deg5_v6.m}{$X_1(34)$} & $21$ & $3$ & $0,0,0,26,32,80$          & $0,0,0,5,4,11$       & $5,7,11$\\
\gitlink{hecke_sieve_deg5_v6.m}{$X_1(35)$} & $25$ & $3$ & $0,0,12,0,36,126$         & $0,0,1,0,3,11$       & $11,13$\\
\gitlink{hecke_sieve_deg5_v6.m}{$X_1(38)$} & $28$ & $3$ & $0,0,3,0,72,87$           & $0,0,1,0,8,11$       & $5,7,11$\\
\gitlink{hecke_sieve_deg5_v6.m}{$X_1(39)$} & $33$ & $2$ & $0,0,0,0,0,12$            & $0,0,0,0,0,1$        & $5,7$\\
\gitlink{hecke_sieve_cached.m}{$X_1(40)$} & $25$ & $3$ & $0,0,0,32,72,80$          & $0,0,0,4,9,10$       & $7,11$\\
\gitlink{hecke_sieve_42_roots.m}{$X_1(42)$} & $25$ & $5$ & $0,0,4,183,648,2286$      & $0,0,2,33,108,392$   & $11,13,17$\\
\gitlink{hecke_sieve_cached.m}{$X_1(49)$} & $69$ & $3$ & $0,0,0,21,63,154$         & $0,0,0,1,3,8$        & $5$\\
\gitlink{hecke_sieve_deg5_v6.m}{$X_1(50)$} & $48$ & $3$ & $0,0,0,0,0,60$            & $0,0,0,0,0,6$        & $7,11$\\
\gitlink{sieve_M2_fast_cached.m}{$X_1(2,22)$} & $16$ & $3$ & $0,0,0,0,60,100$          & $0,0,0,0,12,20$      & $5,7$\\
\gitlink{sieve_M2_tate_cached.m}{$X_1(2,26)$} & $25$ & $3$ & $0,0,0,0,36,122$          & $0,0,0,0,6,26$       & $5,7$\\
\gitlink{sieve_M2_tate_cached.m}{$X_1(2,28)$} & $25$ & $3$ & $0,0,0,0,48,64$           & $0,0,0,0,8,16$       & $5,11$\\
\bottomrule
\end{tabular}
\medskip
\caption{The Hecke sieve computations of \Cref{prop:sieve,prop:49,prop:32,prop:2_26}. The
fourth column lists the number of non-cuspidal closed points of $X_{\F_p}$ of degree
$1,\dots,6$, and the fifth the number of their orbits under the diamond action. Reducible types
such as $[3,3]$ and $[4,2]$ are assembled from these and are processed as well. The sieve leaves
no survivors, except for four divisors on $X_1(32)$ and two on $X_1(2,26)$, which are treated
in \Cref{sec:classes}. In the digital version of the paper, the name of each curve is a link to
the script that runs the sieve on it. For $X_1(33)$, $X_1(42)$, $X_1(49)$ and $X_1(2,22)$ the
computation is split over several scripts, which are listed in the README of the repository.}
\label{table:sieve}
\end{table}

\begin{remark}\label{rem:CC_used}
The cuspidal complements impose the following conditions on $q$. For $(n,p)=(34,3),(38,3)$
the prime $p$ is a primitive root modulo $n$, and for $X_1(35)$, $X_1(39)$ and $X_1(50)$ every
Galois orbit of cusps over $\Q$ stays irreducible modulo $p$, so in these cases every $q$ is
admissible. For $X_1(33)$ at $5$ the only orbit that splits modulo $5$ has degree $10$ and
splits into two places of degree $5$, which cannot occur in a cuspidal complement since
$Y_1(33)(\F_5)=\emptyset$. For $X_1(32)$ at $3$ a cuspidal complement of degree at least $2$ can
contain one of the two places of degree $2$ of an orbit of degree $4$, and for these types we
used only $q=11\equiv3^7\pmod{32}$. For $X_1(40)$ at $3$ the complements of degree at most $2$
impose the condition $q\bmod 40\in\{1,3,9,11,17,19,27,33\}$, and we used only $q=11$ for the
type $[4]$. For $X_1(42)$ at $5$, a cuspidal complement of degree $3$ can be one of the places
of degree $3$ into which two of the orbits of degree $6$ split, and for the type $[3]$ we used
only $q=17\equiv5^5\pmod{42}$. For $X_1(2,22)$ we used only $q=5\equiv 3^8\pmod{44}$ for the
types with a cuspidal complement. For $X_1(2,26)$ and $X_1(2,28)$ the cuspidal complement has
degree at most $1$, and the eighteen (respectively twelve) cusps of degree $1$ over $\F_3$ are
reductions of rational cusps. For the curves $X_1(n)$ the code computes the condition from the
model and checks that the degrees of the reductions of the cusps agree with the prediction of
\cite[Lemma 2.8]{DEHMZ}. For $X_1(2,26)$ and $X_1(2,28)$ it computes the condition from the
description of the cusps by Tate curves, and checks that the degrees of the reductions of the cusps
on the model agree with it. For $X_1(2,22)$ no such computation is needed, as $5\equiv3^8\pmod{44}$.
\end{remark}

\begin{proposition}\label{prop:49}
The group $\torz{49}$ is not a subgroup of the torsion of an elliptic curve over a sextic field.
\end{proposition}

\begin{proof}
We have $\rk J_1(49)(\Q)=0$ by \cite[Theorem 3.1(2)]{DEHMZ}, and Abramovich's bound gives
$\gon_\Q X_1(49)>11.66>6$, so \Cref{prop:main}(a) with $p=3$ excludes the totally cuspidal
reduction. The Galois orbits of cusps over $\Q$ have degrees $1^{21},6^3,21$, and $3$ is a
primitive root modulo $49$, so each of them stays irreducible modulo $3$. By
\Cref{table:sieve} the non-cuspidal part $D$ of the reduction is a single closed point of degree
$4$, $5$ or $6$. In the first two cases the cuspidal complement is supported on the reductions of
rational cusps, and we \gitlink{hecke_sieve_cached.m}{compute} that $[A_5(D)]\neq0$ for the
representatives of the $1$ and $3$ diamond orbits. In the third case there is no cuspidal
complement, and we \gitlink{hecke_sieve_cyclic_commonzero.m}{compute} that $[A_5(D)]\neq0$ also
for the representatives of the $8$ diamond orbits.
\end{proof}

\begin{remark}\label{rem:A2_49}
As an independent check of the case of degree $6$ in the proof of \Cref{prop:49}, we also used
the operator $A_2=T_2-2\dm2-1$. As $q=2$ is not covered by \Cref{prop:tools}(2), this needs a
separate argument, which we give in \Cref{sec:computations}.
\end{remark}

\subsection{The curves over \texorpdfstring{$\Q(\zeta_m)$}{Q(zeta\_m)}}

For $\tg{3}{24}$, $\tg{6}{12}$, $\tg{4}{16}$ and $\tg{4}{20}$ we use the same method over
$\Q(\zeta_m)$, with $d=d_m=3$. We work at a prime $\mathfrak p$ of $\Q(\zeta_m)$ above $p=7,7,3,3$ respectively.
For $m=6$ we have $\dm{13}=1$ on $X_1(6,12)$, so $A_{13}=T_{13}-14$.

\begin{lemma}\label{lem:gon612}
We have $\gon_\C X_1(6,12)>3$.
\end{lemma}

\begin{proof}
The forgetful map $X_1(6,12)\to X_1(3,12)$, $(E,P,Q)\mapsto(E,2P,Q)$, has degree $2$, and the
two curves have genera $9$ and $3$. A map $X_1(6,12)\to\PP^1$ of degree $d\leq3$ cannot factor
through it, as $X_1(3,12)$ is not rational, so the Castelnuovo--Severi inequality \cite[Theorem 3.11.3]{stichtenoth} gives
$9\leq 2\cdot 3+(d-1)\leq 8$, which is a contradiction.
\end{proof}

\begin{proposition}\label{prop:sieve_m}
The groups $\tg{3}{24}$, $\tg{6}{12}$, $\tg{4}{16}$ and $\tg{4}{20}$ are not subgroups of the
torsion of an elliptic curve over a sextic field.
\end{proposition}

\begin{proof}
Let $E$ be an elliptic curve over a sextic field $K$ containing one of the four groups, and let
$X=X_1(m,n)$ be the corresponding curve. In every case $\rk J(\Q(\zeta_m))=0$ by
\cite[Theorem 4.1]{DS}. The indices $[\PSL_2(\Z):\overline{\Gamma_1(n)\cap\Gamma(m)}]=576$, $384$ and $576$
for $X_1(3,24)$, $X_1(4,16)$ and $X_1(4,20)$ give $\gon_\C X>5.71,3.80,5.71$, and $\gon_\C X_1(6,12)>3$ by \Cref{lem:gon612}.
Hence \Cref{prop:main}(a) with $d=3$ shows that $E$ has good reduction at some prime of $K$
above $\mathfrak p$, and the reduction of the trace divisor has a non-cuspidal part $D$.

The residue field $k$ of $\mathfrak p$ is $\F_7$ for $m=3,6$ and $\F_9$ for $m=4$, and the
residue fields of $K$ above $\mathfrak p$ are extensions of $k$ of degree $1$, $2$ or $3$. For $\tg{3}{24}$ and $\tg{6}{12}$ the Hasse
bounds $13$ and $64$ over $\F_7$ and $\F_{49}$ exclude a subgroup of order $72$, for
$\tg{4}{16}$ the Hasse bound over $\F_9$ excludes a subgroup of order $64$, and over $\F_{81}$ a
subgroup of order $64$ forces the Frobenius to be the scalar $9$, so
$E(\F_{81})=E[8]$, which has no point of order $16$. Hence in these three cases $D$ is a single
closed point of degree $3$, with no cuspidal complement. For $\tg{4}{20}$ the closed point can
also have degree $2$, with a cuspidal complement of degree $1$.

The numbers of the non-cuspidal closed points of these degrees and of their diamond orbits are
listed in \Cref{table:sieve_m}, together with the primes $q$ that we use. All of them satisfy
$q\equiv1\pmod m$ and split completely in $\Q(\zeta_m)$, and the standard cusp $c_0$ satisfies
$A_q(c_0)=0$ by \Cref{prop:tools}(3). For $\tg{4}{20}$ and $D$ of degree $2$ we use $q=29$. As
$29\equiv9=3^2\pmod{20}$, the cuspidal complement $C$ satisfies $[A_{29}(C)]=0$ by the discussion
at the beginning of \Cref{sec:sieve}.
For one representative $D$ of each diamond orbit we compute, with the scripts linked in
\Cref{table:sieve_m}, that $[A_q(D)]\neq0$, so no decomposition can occur.
\end{proof}

\begin{table}[H]
\centering
\small
\begin{tabular}{lrrrll}
\toprule
$X$ & $g$ & $k$ & $q$ & non-cuspidal points of degree $1,2,3$ & orbits\\
\midrule
\gitlink{sieve_3_24.m}{$X_1(3,24)$} & $25$ & $\F_7$ & $13$ & $0,0,48$ & $0,0,12$\\
\gitlink{sieve_6_12.m}{$X_1(6,12)$} & $9$  & $\F_7$ & $13$ & $0,0,96$ & $0,0,48$\\
\gitlink{sieve_4_16_p3.m}{$X_1(4,16)$} & $13$ & $\F_9$ & $5$ & $0,0,320$ & $0,0,80$\\
\gitlink{sieve_4_20_p3.m}{$X_1(4,20)$} & $25$ & $\F_9$ & $29,13$ & $0,96,320$ & $0,24,80$\\
\bottomrule
\end{tabular}
\medskip
\caption{The Hecke sieve computations of \Cref{prop:sieve_m}. For $X_1(4,20)$ we use $q=29$
for the points of degree $2$ and $q=13$ for the points of degree $3$. The sieve leaves no
survivors. As in \Cref{table:sieve}, the name of each curve is a link to the script that runs
the sieve on it.}
\label{table:sieve_m}
\end{table}

\begin{remark}
In \Cref{prop:sieve_m} we check the completeness of the enumeration of the closed points by an
independent count of the elliptic curves over the extensions of $k$ of degree $1$, $2$ and $3$
with the required level structure (see \githubbare{independent_m3_counts.m},
\githubbare{independent_m6_counts.m}, \githubbare{independent_m4_p3_counts.m} and
\githubbare{independent_m4_20_p3_counts.m}). For instance, $\#Y_1(4,16)(\F_{9^d})=0,0,960$ for
$d=1,2,3$, which agrees with $320$ closed points of degree $3$. The same check was done for all
the curves of \Cref{table:sieve} (see \Cref{sec:computations}).
\end{remark}

\section{Arguments with cuspidal divisor classes}\label{sec:classes}

The Hecke sieve uses only the fact that $A_q$ kills $J(\Q(\zeta_m))$. When $J(\Q(\zeta_m))$ is known
explicitly, and in particular when it is generated by cuspidal divisor classes, we can instead
compare the classes of the possible reductions with the reduction of $J(\Q(\zeta_m))$ itself. This is
the method of \cite[\S 5.2]{DEHMZ} and \cite[\S 8]{N5}, which we use in the following form.

\begin{lemma}\label{lem:classtest}
Let $X$ be a curve over a number field $F$ with a rational cusp $c_\ast$, let $\mathfrak p$ be a
prime of good reduction of $X$ with residue field $k$, and let $G\subseteq J(k)$ be a subgroup
containing the reduction of $J(F)$. Let $D\in X^{(d)}(F)$, and write its reduction modulo
$\mathfrak p$ as $Z+C$, where $Z$ is effective of degree $e$ with non-cuspidal support and $C$ is
effective and cuspidal. Then
\[
[Z+C-d\,\overline{c_\ast}]\in G.
\]
\end{lemma}

\begin{proof}
This is the reduction of the class $[D-dc_\ast]\in J(F)$.
\end{proof}

So an effective divisor $Z$ of degree $e$ with non-cuspidal support can be the non-cuspidal part
of the reduction of $D$ only if $[Z+C-d\,\overline{c_\ast}]\in G$ for some effective cuspidal divisor
$C$ of degree $d-e$ on $X_k$. To check this for many $Z$ at once, we compute for each $j\geq0$
the set $\mathcal S_j\subseteq J(k)/G$ of the classes of $[j\overline{c_\ast}-C]$, where $C$ runs over
the effective cuspidal divisors of degree $j$ on $X_k$. The condition then says that the image
of $[Z-e\overline{c_\ast}]$ in $J(k)/G$ lies in $\mathcal S_{d-e}$.

In practice $G$ is either the image of the reduction map on $J(F)$, when $J(F)$ is known and
reduction is injective on it, or the subgroup generated by the classes
$[\overline c-\deg(\overline c)\overline{c_\ast}]$ of all closed cusps $\overline c$ of $X_k$, when
$J(F)$ is generated by cuspidal divisor classes. In the latter case $G$ contains the reduction
of $J(F)$ without any assumption on which cusps are rational, and it contains the classes of all
cuspidal divisors of degree $0$ on $X_k$. So every $\mathcal S_j$ is trivial, and the condition is
simply $[Z-e\overline{c_\ast}]\in G$.

Usually we test only one representative $Z$ of each orbit of the diamond operators, and the
numbers of divisors $Z$ that we give are then numbers of orbit representatives. This is
justified by the identity
\[
[\delta Z-e\,\overline{c_\ast}]=\delta[Z-e\,\overline{c_\ast}]+e[\delta\overline{c_\ast}-\overline{c_\ast}]
\qquad(\delta\in\Delta(n)).
\]
In both cases above $G$ is stable under the diamond operators and contains
$[\delta\overline{c_\ast}-\overline{c_\ast}]$, as $\delta c_\ast$ is again an $F$-rational cusp, and
$\delta$ permutes the effective cuspidal divisors of each degree, so the condition holds for $Z$
if and only if it holds for $\delta Z$.

\subsection{The curves \texorpdfstring{$X_1(32)$}{X\_1(32)}, \texorpdfstring{$X_1(2,26)$}{X\_1(2,26)} and \texorpdfstring{$X_1(36)$}{X\_1(36)}}

\begin{proposition}\label{prop:32}
The group $\torz{32}$ is not a subgroup of the torsion of an elliptic curve over a sextic field.
\end{proposition}

\begin{proof}
We have $\rk J_1(32)(\Q)=0$ and $J_1(32)(\Q)=\ClCusp_\Q X_1(32)\simeq
(\torz2)^3\times\torz4\times\torz{120}\times\torz{11640}$ \cite[Theorem 3.1(2), Corollary 4.14 and
Table 2]{DEHMZ}, and $\gon_\Q X_1(32)=8$ \cite[Theorem 4 and Table 1]{DvH}. As in the proof of
\Cref{prop:sieve}, \Cref{prop:main} with $p=3$ excludes the totally cuspidal reduction, and the
Hecke sieve at $p=3$ (see \Cref{table:sieve} and \Cref{rem:CC_used}) excludes every decomposition except
four, three of type $[6]$ and one of type $[3,3]$, none of which has a cuspidal complement. Let
$G$ be the subgroup of $J_1(32)(\F_3)$ generated by the reductions of the rational cuspidal
divisors of degree $0$. We \gitlink{analyze_survivors.m}{compute} that
$G\simeq(\torz2)^3\times\torz4\times\torz{120}\times\torz{11640}$, so $G$ is the reduction of
$J_1(32)(\Q)$, and that none of the four classes $[Z-6\overline{c_\ast}]$ lies in $G$. By
\Cref{lem:classtest} none of them is the reduction of a trace divisor.
\end{proof}

\begin{proposition}\label{prop:2_26}
The group $\tg{2}{26}$ is not a subgroup of the torsion of an elliptic curve over a sextic
field. Hence, by \Cref{lem:4n}, neither is $\torz{52}$.
\end{proposition}

\begin{proof}
Let $X=X_1(2,26)$, a curve of genus $25$. Its Jacobian has rank $0$ over $\Q$ and
$J(\Q)=\ClCusp_\Q X\simeq\torz2\times\torz{14}\times\torz{266}\times\torz{3990}\times\torz{11970}$
\cite[Theorem 3.1(3), Proposition 4.16 and Table 3]{DEHMZ}, and $\gon_\Q X>6$
\cite[proof of Proposition 5.5]{DS}. As in the proof of \Cref{prop:sieve}, the totally cuspidal
reduction is excluded by \Cref{prop:main} with $p=3$, and the Hecke sieve at $p=3$ excludes
every decomposition except two of type $[6]$. The subgroup $G$ of $J(\F_3)$ generated by the
differences of the $18$ rational cusps has invariant factors $[2,14,266,3990,11970]$, so it is
the reduction of $J(\Q)$. We \gitlink{inspect_2_26_classes.m}{compute} that the classes
$[Z-6\overline{c_\ast}]$ of the two remaining divisors do not lie in $G$. As $G$ is stable under
the diamond operators, \Cref{lem:classtest} completes the proof.
\end{proof}

\begin{proposition}\label{prop:36}
The group $\torz{36}$ is not a subgroup of the torsion of an elliptic curve over a sextic field.
\end{proposition}

\begin{proof}
We have $J_1(36)(\Q)=\ClCusp_\Q X_1(36)\simeq\torz{12}\times\torz{252}\times\torz{7812}$
\cite[Theorem 3.1(2), Corollary 4.14 and Table 2]{DEHMZ} and $\gon_\Q X_1(36)=8$
\cite[Theorem 4 and Table 1]{DvH}, so \Cref{prop:main} with $p=5$ excludes the totally cuspidal
reduction. Let $G\subseteq J_1(36)(\F_5)$ be the subgroup generated by the classes
$[\overline c-\deg(\overline c)\overline{c_\ast}]$ of all closed cusps $\overline c$ of
$X_1(36)_{\F_5}$. It contains the reduction of $J_1(36)(\Q)$, and we
\gitlink{class_cuspidal_36_certificate.m}{compute}
\begin{align*}
J_1(36)(\F_5)&\simeq(\torz4)^2\times(\torz{12})^2\times\torz{36}\times(\torz{7812})^2,\\
G&\simeq\torz3\times\torz{36}\times(\torz{7812})^2,
\end{align*}
so $G$ has index $768$. The numbers of non-cuspidal closed points of $X_1(36)_{\F_5}$ of degrees
$1,\dots,6$ are $0,0,48,144,660,2856$, so the non-cuspidal part $Z$ of a possible reduction is a
single point of degree $3,4,5$ or $6$, or the sum of two points of degree $3$. We compute that
$[Z-(\deg Z)\overline{c_\ast}]\notin G$ for all the $48+144+660+2856+1176$ divisors $Z$ of these
types. As $G$ contains the classes of all cuspidal divisors of degree $0$, \Cref{lem:classtest}
shows that none of them occurs.
\end{proof}

\subsection{The sextic points on \texorpdfstring{$X_1(2,24)$}{X\_1(2,24)}}
On $X_1(2,24)$ there are non-cuspidal sextic points, and we determine all of them. The
computation in \Cref{prop:2_24} shows that they form at most two orbits under the diamond
operators, and in \Cref{prop:2_24_points} we show, with an explicit elliptic curve, that this
bound is attained.

\begin{proposition}\label{prop:2_24}
Let $X=X_1(2,24)$. There are $\Q$-rational cuspidal divisors $L_1,L_2$ of degree $6$ on $X$,
given explicitly below, with the following property. If $D$ is an effective $\Q$-rational divisor
of degree $6$ on $X$ with non-cuspidal support, then $\delta D\sim L_i$ for some $\delta\in\Delta(24)$
and some $i\in\{1,2\}$. Each of the linear systems $|L_1|$ and $|L_2|$ contains at most one
effective divisor.
\end{proposition}

\begin{proof}
The curve $X$ has genus $13$, its Jacobian has rank $0$ and
$J(\Q)=\ClCusp_\Q X\simeq\torz2\times\torz4\times\torz4\times\torz{120}\times\torz{240}$
\cite[Theorem 3.1(3), Proposition 4.16 and Table 3]{DEHMZ}, and $\gon_{\F_5}X>6$
\cite[proof of Proposition 5.5]{DS}, so $\gon_\Q X>6$. We work at $p=5$.

Let $D$ be an effective $\Q$-rational divisor of degree $6$ on $X$ with non-cuspidal support. If
$D$ were not a single closed point of degree $6$, its support would contain a closed point of
degree at most $3$, which would give an elliptic curve over a number field of degree at most $3$
with a subgroup isomorphic to $\tg{2}{24}$, which is impossible
\cite{mazur,kamienny,DEHMZ}. So $D$ is a closed point of degree $6$, and hence the trace
divisor of a point of $Y_1(2,24)(K)$, where $K$ is its residue field. In particular the totally
cuspidal reduction of $D$ is excluded by \Cref{prop:main}.

We use the normalized universal curve $y^2=x^3+cx^2+dx$ with $P=(0,0)$ of order $2$ and
$Q=(b,b)$ of order $24$, where $d=b(1-b-c)$. The functions $b$, $d$, $c^2-4d$ and $x(jQ)$ for
$1\leq j\leq 12$ are units on $Y_1(2,24)_{\Z[1/6]}$, since $jQ$ for $j\leq12$ is never equal to
$O$ or $P$ in any fibre. It follows that their orders of vanishing at the cusps do not change under
reduction modulo $5$, and grouping the cusps by these orders of vanishing gives $19$ blocks
$B_1,\dots,B_{19}$, each of which is a sum of Galois orbits of cusps over $\Q$, of degrees
$1,1,1,2,2,1,4,4,4,4,2,2,1,4,2,2,1,1,1$, where we order the blocks increasingly in the
lexicographic order of their vectors of orders of vanishing
$(\ord b,\ord d,\ord(c^2-4d),\ord x(Q),\dots,\ord x(12Q))$. The block $B_1$ is a single rational
cusp, at which $b$ has order $-22$. In this proof we take $c_\ast:=B_1$ in \Cref{lem:classtest}, and we put
\begin{align*}
L_1&=15B_1-B_7-B_8+B_{11}-B_{14}+B_{17},\\
L_2&=3B_1+B_3+B_4+B_{13}-B_{16}+B_{19}.
\end{align*}
Let $G\subseteq J(\F_5)$ be the subgroup generated by the classes
$[\overline{B_j}-\deg(B_j)\overline{c_\ast}]$. We \gitlink{fine_2_24_cusp_reduction.m}{compute} that
$G$ has invariant factors $[2,4,4,120,240]$, so $G$ is the reduction of $J(\Q)$, and that
$J(\F_5)/G\simeq(\torz2)^2\times\torz{24}\times\torz{48}$. The sets $\mathcal S_j$ defined after
\Cref{lem:classtest} have $1,5,23,75$ elements for $j=0,1,2,3$.

Over $\F_5$ the numbers of non-cuspidal closed points of degrees $1,\dots,6$ are
$0$, $0$, $16$, $96$, $672$, $2712$, so the non-cuspidal part $Z$ of the reduction of $D$ has type
$[3]$, $[4]$, $[5]$, $[6]$ or $[3,3]$. Up to the diamond operators there are $4$, $32$, $168$,
$684$ and $64$ such $Z$, respectively, and we compute that the condition
$[Z-(\deg Z)\overline{c_\ast}]\bmod G\in\mathcal S_{6-\deg Z}$ holds for exactly two of them, both of type $[6]$, whose classes are those
of $L_1$ and $L_2$.

So, after applying a diamond operator, $D$ reduces to one of these two points, and then $[D-L_i]$ is a
rational torsion class with trivial reduction, so $D\sim L_i$. As $\gon_\Q X>6$, each of the
linear systems $|L_1|$ and $|L_2|$ contains at most one effective divisor.
\end{proof}

\subsection{The group \texorpdfstring{$\tg{3}{18}$}{Z/3Z x Z/18Z}}
For $\tg{3}{18}$ the Hecke sieve at a single prime does not suffice, and we do not know
$J_1(3,18)(\Q(\zeta_3))$ explicitly. Instead we compare the reductions at two primes.

\begin{lemma}\label{lem:matched}
Let $X$ be a curve over a number field $F$ with $J(F)$ finite and with a rational cusp $c_\ast$.
Let $\mathfrak p_1,\mathfrak p_2$ be primes of good reduction, with residue fields $k_1,k_2$, such
that reduction modulo $\mathfrak p_1$ is injective on $J(F)$, let $H_i\subseteq J(k_i)$ be the
image of $\ClCusp_F X$ under reduction modulo $\mathfrak p_i$, and let $\psi:H_1\to H_2$ be the
composite of the inverse of the reduction map $\ClCusp_F X\to H_1$ with the reduction map
$\ClCusp_F X\to H_2$. If $J(k_1)/H_1$ has exponent $e$, then for every $D\in X^{(d)}(F)$, with
reductions $\overline D_i$ and $\overline{c_{\ast,i}}$ of $D$ and $c_\ast$ modulo $\mathfrak p_i$, we have
\[
\psi\bigl(e[\overline D_1-d\overline{c_{\ast,1}}]\bigr)=e[\overline D_2-d\overline{c_{\ast,2}}].
\]
\end{lemma}

\begin{proof}
Let $y=[D-dc_\ast]\in J(F)$. The reduction of $ey$ modulo $\mathfrak p_1$ lies in $H_1$, so it is the
reduction of some $c\in\ClCusp_F X$. Then $ey-c$ has trivial reduction modulo $\mathfrak p_1$, so
$ey=c$. Reducing modulo $\mathfrak p_2$ gives the claim.
\end{proof}

\begin{proposition}\label{prop:3_18}
The group $\tg{3}{18}$ is not a subgroup of the torsion of an elliptic curve over a sextic field.
\end{proposition}

\begin{proof}
Let $X=X_1(3,18)$ over $\Q(\zeta_3)$, a curve of genus $10$. We have $\rk J(\Q(\zeta_3))=0$
\cite[Theorem 4.1]{DS}, and the index $[\PSL_2(\Z):\overline{\Gamma_1(18)\cap\Gamma(3)}]=324$ gives
$\gon_\C X>3.21>3$. We work at primes
$\mathfrak p_1$ and $\mathfrak p_2$ of $\Q(\zeta_3)$ above $7$ and $13$, with residue fields $\F_7$ and
$\F_{13}$. By \Cref{prop:main}, applied to $\mathfrak p_1$ and to $\mathfrak p_2$, the reduction
of the trace divisor $D$ modulo $\mathfrak p_i$ has a non-cuspidal part. The numbers of non-cuspidal closed points of degrees $1,2,3$ are
$0,27,144$ over $\F_7$ and $0,54,846$ over $\F_{13}$, so this non-cuspidal part is either a
single point of degree $3$ or a point of degree $2$, in which case the cuspidal complement is a
cusp of degree $1$, which is the reduction of a $\Q(\zeta_3)$-rational cusp. The diamond operators act
on these points with orbits of length $3$. At $\mathfrak p_1$ we
\gitlink{diagonal_3_18_cusp_classes_certificate.m}{compute} that $A_{13}$ excludes all $9$ orbits of points of degree
$2$ and $36$ of the $48$ orbits of points of degree $3$, and at $\mathfrak p_2$ that $A_7$
excludes all $18$ orbits of points of degree $2$ and $270$ of the $282$ orbits of points of
degree $3$. This leaves $36$ possible reductions at each prime.

To compare the two primes we need to identify the reductions of the cuspidal divisor classes.
We use the model of $X$ of \cite{DSmodels}, as implemented in \cite{mdmagma}, on which the
universal curve is $y^2+a_1xy+a_3y=x^3$ with $P=(0,0)$ of order $3$, normalized by
$y(Q)=x(Q)^2$ for the point $Q$ of order $18$. A point of this curve with $x=0$ is $\pm P$, and
one with $y=0$ is $P$, so $x(Q)$ and $y(Q)$ have no zeros or poles on $Y_1(3,18)$, and the
normalization makes $a_1$ and $a_3$ regular functions on $Y_1(3,18)$ over
$\Z[\zeta_3,\tfrac16]$. The functions
\[
a_3,\quad a_1^3-27a_3,\quad x(kQ+jP),\ y(kQ+jP)\ \ (1\leq k\leq 9,\ 0\leq j\leq 2),
\]
\[
x(kQ)-x(kQ+P)\ \ (1\leq k\leq 9),\qquad x(kQ+P)-x(kQ+2P)\ \ (1\leq k\leq 8)
\]
are units on $Y_1(3,18)$ over $\Z[\zeta_3,\tfrac16]$. Indeed, the discriminant of the curve is
$a_3^3(a_1^3-27a_3)$, the points $kQ+jP$ are never equal to $O$ or $\pm P$, and in every fibre
$kQ+P\neq\pm kQ$ and $kQ+2P\neq\pm(kQ+P)$ for these $k$. So grouping the cusps by their orders
of vanishing at these functions gives
the same $24$ blocks, of degrees $1^{18},3^6$, at both primes, each of them the reduction of a
Galois orbit of cusps over $\Q(\zeta_3)$. We compute that the images $H_1,H_2$ of $\ClCusp_{\Q(\zeta_3)}X$ both have
invariant factors $[3,3,3,3,63,63,63]$ and that $J(\F_7)/H_1\simeq\torz3\times\torz3\times\torz9$.
As $7$ is odd, unramified in $\Q(\zeta_3)$ and a prime of good reduction for $X$, reduction modulo
$\mathfrak p_1$ is injective on $J(\Q(\zeta_3))=J(\Q(\zeta_3))_\tors$. Hence the map $\psi$ of
\Cref{lem:matched} sends the reduction modulo $\mathfrak p_1$ of the class of a combination of the
blocks of degree $0$ to its reduction modulo $\mathfrak p_2$. By \Cref{lem:matched} with $e=9$, a trace divisor would give a pair of possible reductions
$(\overline D_1,\overline D_2)$ with $\psi(9[\overline D_1-3\overline{c_{\ast,1}}])=
9[\overline D_2-3\overline{c_{\ast,2}}]$, and we compute that none of the $36\cdot36$ pairs satisfies
this.
\end{proof}

\subsection{The group \texorpdfstring{$\tg{2}{30}$}{Z/2Z x Z/30Z}}
For $X=X_1(2,30)$ we push the divisors forward to quotient curves. Let $Y=X_1(30)$ and consider
the three morphisms
\[
\pi_1(E,P,Q)=(E,Q),\qquad \pi_2(E,P,Q)=(E,P+16Q),\qquad \pi_3(E,P,Q)=(E,P+Q)
\]
from $X$ to $Y$. They are defined over $\Q$, they map cusps to cusps and they commute with the
diamond operators. This is because $15Q$, $P$ and $P+15Q$ are the non-zero points of $E[2]$, and adding
$16Q$, which has order $15$, to them gives points of order $30$. These three maps do not
suffice, and we also use the quotient $W$ of $X$ by the automorphism $(E,P,Q)\mapsto(E,15Q,P+Q)$
of order $3$, which permutes the three non-zero points of $E[2]$ cyclically. The curve $W$ sees a
part of $J(X)$ that is invisible on $Y$ (see \Cref{rem:2_30}). It is the double cover of
$X_1(15)$ given by a square root of the discriminant of the
universal elliptic curve over $X_1(15)$ in Tate normal form, every cusp of $X_1(15)$ is totally
ramified in it, and it has genus $9$, as we \gitlink{a3_m2_30_model.m}{compute}.

\begin{proposition}\label{prop:2_30}
The group $\tg{2}{30}$ is not a subgroup of the torsion of an elliptic curve over a sextic field.
Hence, by \Cref{lem:4n}, neither is $\torz{60}$.
\end{proposition}

\begin{proof}
We have $\rk J(\Q)=0$ \cite[Theorem 4.1]{DS}, and $\gon_\Q X>6$ by
\cite[Lemma 5.6 and Proposition 5.7]{DS}, so \Cref{prop:main} with $p=7$ excludes the totally cuspidal
reduction. Write the reduction modulo $7$ of the trace divisor as $Z'+C$, with $Z'$ non-cuspidal
of degree $e>0$.

By \cite[Theorem 3.1(2), Corollary 4.14 and Table 2]{DEHMZ},
$J(Y)(\Q)=\ClCusp_\Q Y\simeq\torz4\times\torz{8160}$, and we
\gitlink{m2_30_quotient_prefilter.m}{compute} that $J(Y)(\F_7)\simeq(\torz4)^3\times\torz{136}\times\torz{16320}$ and that the reduction $G$ of
$J(Y)(\Q)$ has index $4352$. Applying \Cref{lem:classtest} on $Y$ to the push-forwards
$\pi_{i*}$ of the trace divisor, with every cuspidal complement on $Y$ allowed, gives the
condition
\[
[\pi_{i*}Z'-e\overline{c_\ast}]\bmod G\in\mathcal S_{6-e}\qquad(i=1,2,3),
\]
where the sets $\mathcal S_0,\dots,\mathcal S_4$ have $1,9,41,67,68$ elements. The numbers of
non-cuspidal closed points of $X_{\F_7}$ of degrees $1,\dots,6$ are $0,24,192,732,2736,20000$.
We \gitlink{m2_30_quotient_prefilter.m}{compute} that of the $22706$ possible non-cuspidal parts $Z'$, up
to the diamond operators, $68$ satisfy the condition for $\pi_1$, and we
\gitlink{m2_30_other_quotients.m}{check} that $10$ of these satisfy it for all three maps. These $10$ divisors have degree $6$. Four are single points of degree $6$ and six
are sums of three points of degree $2$.

To exclude them we use the quotient $W$. By \cite[Proposition 4.16]{DEHMZ}, every element of
$J(\Q)$ is the class of a (not necessarily effective) $\Q$-rational cuspidal divisor of degree
$0$. Push-forward maps such divisors on $X$ to $\Q$-rational cuspidal divisors of degree $0$ on
$W$, so the push-forward to $W$ of the class of any rational divisor of degree $0$ on $X$ is the
class of a rational cuspidal divisor on $W$. Let $H\subseteq J(W)(\F_7)$ be the reduction of the group generated by the rational
cuspidal divisor classes of degree $0$ on $W$. We \gitlink{filter_a3_m2_30.m}{compute}
\[
J(W)(\F_7)\simeq(\torz2)^6\times\torz{102}\times\torz{10608},\qquad H\simeq(\torz2)^5\times\torz8,
\]
and that the classes of the push-forwards to $W$ of the $10$ remaining divisors, based at
$6\overline{c_\ast}$ for a rational cusp $c_\ast$ of $W$, do not lie in $H$. This excludes them.
\end{proof}

\begin{remark}\label{rem:2_30}
The covering $X\to X_1(15)$, $(E,P,Q)\mapsto(E,2Q)$, is Galois with group $\SL_2(\F_2)\simeq S_3$,
acting on the ordered bases $(P,15Q)$ of $E[2]$, and $W$ is its quotient by the subgroup $A_3$
of order $3$. The space $H^0(X,\Omega^1)$ decomposes into $1$ copy of the trivial representation
of $S_3$, $8$ copies of the sign representation and $8$ copies of the $2$-dimensional one, as the
genera $1$, $9$, $9$ and $25$ of $X_1(15)$, $W$, $Y$ and $X$ show. Each $\pi_i$ is invariant under
a subgroup of order $2$ of $S_3$, so the differentials pulled back from $Y$ involve only the
trivial and the $2$-dimensional representations, while those on $W$ are the invariants of $A_3$,
which are the trivial and the sign representations. So the test on $W$ detects the part of
$J(X)$, of dimension $8$, on which $S_3$ acts by the sign character and which is invisible on
$Y$.
\end{remark}

\subsection{The group \texorpdfstring{$\torz{54}$}{Z/54Z}}
Let $X=X_1(54)$ and $Y=X_1(27)$, and let $\pi:X\to Y$ be the map $(E,P)\mapsto(E,2P)$. It is
defined over $\Q$, it has degree $3$, and it maps cusps to cusps and non-cuspidal points to
non-cuspidal points. Let $z$ be a non-cuspidal closed point of $Y_{\F_5}$, with elliptic curve
$E_z$ over the residue field of $z$. The closed points of $X_{\F_5}$ above $z$ correspond to the
Galois orbits of points of order $2$ of $E_z$, so such a point $x$ has degree $\deg z$ or
$2\deg z$ if $E_z$ has a point of order $2$ over the residue field of $z$, and degree $3\deg z$
otherwise, and $\pi_*x=(\deg x/\deg z)\,z$.

\begin{proposition}\label{prop:54}
The group $\torz{54}$ is not a subgroup of the torsion of an elliptic curve over a sextic field.
\end{proposition}

\begin{proof}
Suppose that $E/K$ has a point $P$ of order $54$. We have $\rk J_1(54)(\Q)=0$
\cite[Theorem 3.1(2)]{DEHMZ}, and the index $[\PSL_2(\Z):\overline{\Gamma_1(54)}]=972$ gives
$\gon_\Q X\geq\gon_\C X>9.64>6$, so \Cref{prop:main} with $p=5$ excludes the totally cuspidal
reduction. Write the reduction modulo $5$ of the trace divisor $D$ of $(E,P)$ as $Z'+C$, with
$Z'$ non-cuspidal of degree $e>0$.

By \cite[Corollary 4.14 and Table 2]{DEHMZ},
$J_1(27)(\Q)=\ClCusp_\Q Y\simeq\torz3\times\torz3\times\torz{52497}$. Let
$G\subseteq J_1(27)(\F_5)$ be the subgroup generated by the classes
$[\overline c-\deg(\overline c)\overline{c_\ast}]$ of all closed cusps $\overline c$ of
$Y_{\F_5}$. It contains the reduction of $J_1(27)(\Q)$, and we
\gitlink{certify_54_mod5.m}{compute}
\[
J_1(27)(\F_5)\simeq\torz3\times\torz3\times\torz{342175446},\qquad
G\simeq\torz3\times\torz3\times\torz{52497},
\]
so $G$ has index $6518$. As $\pi$ is defined over $\Q$, the divisor $\pi_*D$ on $Y$ is effective,
$\Q$-rational and of degree $6$, and its reduction modulo $5$ is $\pi_*Z'+\pi_*C$, where $\pi_*Z'$
is non-cuspidal of degree $e$ and $\pi_*C$ is cuspidal. So \Cref{lem:classtest} on $Y$ gives
$[\pi_*Z'+\pi_*C-6\overline{c_\ast}]\in G$, and as $G$ contains the classes of all cuspidal
divisors of degree $0$, it follows that $[\pi_*Z'-e\overline{c_\ast}]\in G$.

We also compute that the numbers of non-cuspidal closed points of $Y_{\F_5}$ of degrees
$1,\dots,6$ are $0,9,57,171,612,2700$, and that the elliptic curves of
$0,0,30,171,396,1890$ of them have a point of order $2$ over the residue field. As $e\leq6$, it
follows that $\pi_*Z'$ is either $3z$ with $z$ of degree $2$, or a single point $z$ of degree
$3,4,5$ or $6$, or $z+z'$ with $z$ and $z'$ of degree $3$ (possibly equal), where in the last
two cases the elliptic curve of each point has a point of order $2$ over its residue field. We
\gitlink{certify_54_mod5.m}{compute} that $[\pi_*Z'-e\overline{c_\ast}]\in G$ for only six of the
$9+30+171+396+1890+465$ divisors of these types, namely $z_i+z_j$ with $1\leq i\leq j\leq3$, where
$z_1,z_2,z_3$ are the points of degree $3$ of $Y_{\F_5}$ whose elliptic curves have
$j$-invariant $2$. For a prime $\wp$ of $K$ above $5$, the reduction of $(E,P)$ modulo $\wp$ is a point of $X$ with
values in $\F_\wp$, so it lies on a closed point of $X_{\F_5}$, which is non-cuspidal if and only if
$E$ has good reduction at $\wp$. The reduction of $D$ modulo $5$ is a sum of positive multiples of
these closed points, so the support of $Z'$ consists of those for the primes $\wp$ at which $E$
has good reduction, and as $e>0$ there is at least one such prime. As $\pi$ maps each of them to
one of $z_1,z_2,z_3$, the reduction of $E$ modulo each such $\wp$ has $j$-invariant $2$.

Let $\wp$ be such a prime. The curve $E'=E/\dm{27P}$ also has a point of order $54$
over $K$: if $\phi:E\to E'$ is the quotient isogeny, then $\phi(2P)$ has order $27$, and
$\ker\hat\phi=\phi(E[2])$ is generated by a $K$-rational point of order $2$. So the above
applies to $E'$ as well, and as $E'$ has good reduction at $\wp$, its reduction modulo $\wp$
also has $j$-invariant $2$. But this reduction is the quotient of the reduction of $E$ by a
subgroup of order $2$. The curve $y^2=x^3+x+1$ has $j$-invariant $1728\cdot4/31=2$, so over
$\overline\F_5$ the reduction of $E$ is isomorphic to it. Its points of order $2$ are $(r,0)$ with
$r^3+r+1=0$, and as $x^3+x+1$ has no root in $\F_5$, $r$ has degree $3$ over $\F_5$. Moving
$(r,0)$ to $(0,0)$ gives the model $y^2=x(x^2+3rx+3r^2+1)$, and the quotient of a curve
$y^2=x(x^2+ax+b)$ by $\dm{(0,0)}$ is $y^2=x(x^2-2ax+a^2-4b)$. So the quotient by
$\dm{(r,0)}$ is $y^2=x(x^2+4rx+2r^2+1)$, whose $j$-invariant is
\[
\frac{1}{r^2+r+1}=2r^2+r+4\notin\F_5,
\]
which is a contradiction.
\end{proof}

\section{Formal immersion at the prime \texorpdfstring{$2$}{2}}\label{sec:formal}

For $n=51$, $63$ and $65$ every elliptic curve over a sextic field with a point of order $n$ has
bad reduction at all primes above $2$ (for $n=65$ this is proved in \Cref{sec:57_65}). But
\Cref{prop:main} does not apply at $p=2$: for $n=51$ Abramovich's bound only gives
$\gon_\C X_1(51)>11.42$, which does not show that the gonality is larger than $12$, and for
$n=63,65$ there are newforms of level dividing $n$ with positive analytic rank, so we cannot use
that $J_1(n)(\Q)$ is finite. Neither does \Cref{prop:several}, which with
$\psi=\dm a-1$ would need a lower bound larger than $24$, while Abramovich's bound only gives
$\gon_\C X_1(n)>11.42,17.13,19.99$. We exclude the totally cuspidal reduction at
$2$ by a formal immersion argument.

The difficulty at $p=2$ is that reduction is in general not injective on the torsion of the
Jacobian. The basic example is the group $\mu_2=\{\pm1\}$ over $\Z_2$: the point $-1$ has order $2$,
but it reduces to $1$ modulo $2$. The following lemma shows that a suitable Hecke operator kills
such points.

\begin{lemma}\label{lem:dyadic}
Let $A$ be an abelian scheme over $\Z_2$, and let $T$ and $\delta$ be endomorphisms of $A$ such
that $T=F+\delta V$ on the special fibre $A_{\F_2}$, where $F$ and $V$ are the Frobenius and the
Verschiebung of $A_{\F_2}$. If $z\in A(\Q_2)$ satisfies $2z=0$ and reduces to $0$ modulo $2$, then
$(T-\delta)z=0$.
\end{lemma}

\begin{proof}
Choose coordinates $X=(X_1,\dots,X_g)$ of the formal group of $A$, that is, functions vanishing
on the zero section whose differentials at the origin of $A_{\F_2}$ form a basis of the cotangent
space. Then the points of $A(\Q_2)$ that reduce to $0$ are the points with coordinates in $2\Z_2$
(see the proof of Lemma$'$ in \cite[Appendix]{katz}). Write
\[
[2](X)=2X+Q(X)+R(X),
\]
where the coordinates of $Q$ are homogeneous of degree $2$ and $R$ has only terms of degree at
least $3$, all with coefficients in $\Z_2$. Let $M$ be a matrix over $\Z_2$ whose reduction
$\overline M$ is the matrix of the linear part of $V$. As $[2]=V\circ F$ on $A_{\F_2}$ and $F$
squares the coordinates, the reduction of $Q$ is $\overline Q(X)=\overline M(X_1^2,\dots,X_g^2)^t$.

Write the coordinates of $z$ as $2u$ with $u\in\Z_2^g$. As $2z=0$, we have
\[
0=[2](2u)\equiv4u+4Q(u)\pmod 8,
\]
so $\overline u+\overline M(\overline u_1^2,\dots,\overline u_g^2)^t=0$, and as $a^2=a$ for
$a\in\F_2$, we get $(\overline M-1)\overline u=0$. Put $\alpha=T-\delta$. As the linear part of $F$
is $0$, the linear part of $\alpha$ on $A_{\F_2}$ is $\overline N(\overline M-1)$, where
$\overline N$ is the linear part of $\delta$. Hence the linear term of $\alpha$ at $2u$ is
divisible by $4$, and so are all the other terms, so the coordinates of $\alpha z$ lie in $4\Z_2$.

Finally, a non-zero point $v$ with coordinates in $4\Z_2$ is not killed by $2$. Indeed, if
$s\geq2$ is the least valuation of a coordinate of $v$, then in such a coordinate the linear
term of $[2](v)$ has valuation $s+1$, and all the other terms have valuation at least $2s>s+1$.
As $2\alpha z=\alpha(2z)=0$, we get $\alpha z=0$.
\end{proof}

We apply this to quotients of $X_1(n)$. Let $n$ be odd and let $H\subseteq\Delta(n)$ be a
subgroup that acts freely on the geometric points of $X_1(n)$ in characteristic $0$ and in
characteristic $2$, cusps included. Then $Y:=X_1(n)/H$ is a smooth curve over $\Z_{(2)}$ with good
reduction at $2$, the Hecke correspondence $T_2$ and the diamond operators descend to $Y$, and
the Eichler--Shimura relation $T_2=F+\dm2V$ holds on $J(Y)_{\F_2}$. We write $c_0$ for the image
of the standard cusp, and $\Cart$ for the Cartier operator on the differentials of $Y_{\F_2}$:
if $x$ is a separating variable, every differential can be written uniquely as
$(f^2+g^2x)\,dx$ with $f,g$ rational functions on $Y_{\F_2}$, and
$\Cart\bigl((f^2+g^2x)\,dx\bigr)=g\,dx$. Equivalently, in terms of a local parameter $t$ we have
$\Cart\bigl(\sum_na_nt^n\,dt\bigr)=\sum_na_{2n+1}^{1/2}t^n\,dt$. The operator $\Cart$ preserves
$H^0(Y_{\F_2},\Omega^1)$ and is $\F_2$-linear on it, and at a cusp at which $q$ is a local
parameter it acts on $q$-expansions by
$\sum_na_nq^n\,dq/q\mapsto\sum_na_{2n}^{1/2}q^n\,dq/q$. If all $a_n$ lie in $\F_2$, this is
Atkin's operator $U_2\colon\sum_na_nq^n\mapsto\sum_na_{2n}q^n$ applied to $\sum_na_nq^n$.

Every endomorphism $\psi$ of $J(Y)$ extends to the N\'eron model of $J(Y)$ over $\Z_2$, which is
an abelian scheme, and we write $\psi^*$ for its action on the cotangent space
$H^0(Y_{\F_2},\Omega^1)$ of the special fibre. As the Frobenius pulls back differentials to $0$
and the Verschiebung pulls back to $\Cart$ \cite{oda}, the Eichler--Shimura relation gives
\begin{equation}\label{eq:T2star}
(T_2-\dm2)^*=(\Cart-1)\dm2^*\quad\text{on } H^0(Y_{\F_2},\Omega^1).
\end{equation}
The following criterion combines three ingredients. The elements $\psi_i$, which have finite
image on $J(Y)(\Q)$, turn the difference of a rational divisor and a cuspidal divisor with the
same reduction into a torsion class. \Cref{lem:dyadic} shows that $T_2-\dm2$ kills this torsion
class. Finally, by \eqref{eq:T2star}, the image of the cotangent map at $\overline C$ of the
resulting morphism $Y^{(d)}\to J(Y)^r$ is the restriction to $\overline C$ of the space $U$ below,
so the condition on $U$ says that this morphism is a formal immersion at $\overline C$. Since the
target is a power of $J(Y)$ itself rather than a quotient of it, this cotangent map can be
computed on $H^0(Y_{\F_2},\Omega^1)$.

\begin{proposition}\label{prop:cartier}
Let $Y=X_1(n)/H$ be as above, let $\psi_1,\dots,\psi_r\in\Z[\Delta(n)/H]$ be such that each
$\psi_i\bigl(J(Y)(\Q)\bigr)$ is finite, and put
\[
U:=(\Cart-1)\Bigl(\sum_{i=1}^r\psi_i^*H^0(Y_{\F_2},\Omega^1)\Bigr)\subseteq
H^0(Y_{\F_2},\Omega^1).
\]
Let $\overline C$ be an effective cuspidal divisor of degree $d$ on $Y_{\F_2}$ such that the
restriction map $U\to H^0(\overline C,\Omega^1|_{\overline C})$ is surjective. Then no effective
$\Q$-rational divisor $D$ of degree $d$ on $Y$ with non-cuspidal support reduces to $\overline C$
modulo $2$.
\end{proposition}

\begin{proof}
Suppose that $D$ reduces to $\overline C$. The cuspidal subscheme of $Y$ is finite \'etale at $2$, so
$\overline C$ lifts uniquely to a cuspidal divisor $C$, which is fixed by a Frobenius element at
$2$, so it is defined over the decomposition field $L\subseteq\Q(\zeta_n)$ of a prime
$\lambda$ above $2$, which satisfies $L_\lambda=\Q_2$. For each $i$ the class
$\psi_i[D-C]=\psi_i[D-dc_0]-\psi_i[C-dc_0]$ is torsion, by the assumption on $\psi_i$ and by
Manin--Drinfeld, and it reduces to $0$ modulo $\lambda$. By Parent's result it is killed by $2$,
so by \Cref{lem:dyadic}, applied to the N\'eron model of $J(Y)$ over $\Z_2$ with $T=T_2$ and
$\delta=\dm2$, we get $(T_2-\dm2)\psi_i[D-C]=0$.

Consider the morphism $f:Y^{(d)}\to J(Y)^r$ over $\Z_2$ given by
$f(B)=\bigl((T_2-\dm2)\psi_i[B-dc_0]\bigr)_i$. Then $f(D)=f(C)$, and $D$ and $C$ have the same
reduction $\overline C$. By \eqref{eq:T2star}, and since $\dm2^*$ is invertible and commutes with $\Cart$ and with the $\psi_i^*$, the image of the cotangent map of $f$ at
$\overline C$ is the restriction of $U$ to $\overline C$, as the cotangent map of
$Y^{(d)}\to J(Y)$, $B\mapsto[B-dc_0]$, at $\overline C$ is the restriction map
$H^0(Y_{\F_2},\Omega^1)\to H^0(\overline C,\Omega^1|_{\overline C})$. By assumption it is
surjective, so $f$ is a formal immersion at $\overline C$, and by \cite[Lemma 4.2]{DKSS} (which is
due to Parent \cite[Lemme 4.13]{parentformal}) $f$ is injective on the residue class of
$\overline C$. Hence $D=C$, which is a contradiction since $D$ is non-cuspidal.
\end{proof}

As in \Cref{sec:rank0}, we check the surjectivity of the restriction maps by computing, for each
closed cusp $P$ of $Y_{\F_2}$ and each $j$, the subspace $\Omega(jP)$ of $H^0(Y_{\F_2},\Omega^1)$,
and then, for each $\overline C$, the rank $\dim U-\dim(U\cap\Omega(\overline C))$ of the
restriction map on $U$, whose kernel is $U\cap\Omega(\overline C)$. Note that \Cref{prop:cartier} uses the rank of $J(Y)$ only through the
operators $\psi_i$, and that we can take $r=1$ and $\psi_1=1$ when $J(Y)(\Q)$ is finite.

\begin{remark}\label{rem:kamienny}
The condition in \Cref{prop:cartier} is the usual formal immersion criterion, written in terms of
cotangent spaces. When $\overline C$ is a sum of rational cusps it corresponds, by the duality
between the Hecke algebra and $q$-expansions, to the linear independence of suitable Hecke
operators modulo $2$, as in the criterion of Kamienny \cite{kamiennyHD} for $X_0(p)$ and its
generalization by Parent \cite[Th\'eor\`eme 1.10 and Proposition 2.7]{parent} (see
\cite[Proposition 5.3]{DKSS}), and when $\overline C$ is a sum of distinct points it is
\cite[Lemma 8.3]{DKSS}. The usual way of dealing with the torsion in the kernel of reduction at
$2$ is Parent's trick \cite[Corollary 5.2]{DKSS}, which multiplies the Hecke operator by an
operator at an odd prime $q$ that kills all the rational torsion. This could also be used here,
with $A_q$ for an odd prime $q$ that splits completely in the field $L$ of the proof of
\Cref{prop:cartier}, that is, with $q$ congruent modulo $n$ to a power of $2$, by the argument of
\Cref{prop:tools}(2). We use $T_2-\dm2$ instead for computational reasons. Its action on
differentials modulo $2$ is given by the Cartier operator, which is computed directly from a
model of $Y_{\F_2}$, whereas $A_q$ requires the Hecke correspondence $T_q$ on $Y$ for a prime $q$
that need not be small (for $n=63$ the smallest such $q$ is $67$). Moreover, the condition of
\Cref{prop:cartier} is checked in the same way at every effective cuspidal divisor
$\overline C$, including those supported on cusps that are not rational and those with
multiplicities.
\end{remark}

\begin{lemma}\label{lem:quotients_51_63}
Let $Y_{51}:=X_1(51)/\dm{16}$ and $Y_{63}:=X_1(63)/\dm8$. The diamond operators $\dm{16}$ and
$\dm8$ have order $2$ and act freely on $X_1(51)$ and $X_1(63)$, also in characteristic $2$, so
$Y_{51}$ and $Y_{63}$ have good reduction at $2$ and have genus $33$ and $49$. The function field
of $Y_{51}$ is obtained by adjoining to the function field of $X_1(17)$ a root of the
$3$-division polynomial of the universal elliptic curve, and that of $Y_{63}$ by adjoining to the
function field of $X_1(9)$ a root of the $7$-division polynomial of the universal elliptic curve.
\end{lemma}

\begin{proof}
Write $(r,s,a)=(17,3,16)$ and $(9,7,8)$ in the two cases, so that $a\equiv-1\pmod r$ and
$a\equiv1\pmod s$. A point of order $rs$ corresponds to a pair of points $(P_r,P_s)$ of orders
$r$ and $s$, and modulo the common sign, $\dm a$ identifies $(P_r,P_s)$ with $(P_r,-P_s)$. So
$Y$ parametrizes pairs $(E,P_r)$ together with the $x$-coordinate of $P_s$, which gives the
description of the function field. If $\dm a$ fixed a non-cuspidal geometric point, there would
be an automorphism $u$ of an elliptic curve with $u(P)=aP$. As $a^2=1$, the automorphism $u^2$
fixes $P$, which has order at least $5$, so $u^2=1$, as a non-trivial automorphism of an
elliptic curve fixes no point of order at least $5$. Hence $u=\pm1$, which acts on $P$ as
$\pm1\neq a$. By
\Cref{sec:prelim}, a cusp fixed by $\dm a$ would require $a$ to be congruent to the same sign
modulo every prime dividing $n$, while $a\equiv-1\pmod r$ and $a\equiv1\pmod s$, so $\dm a$
fixes no cusp. As the cuspidal subscheme is finite \'etale over $\Z[\zeta_n,\tfrac1n]$ and $n$
is odd, $\dm a$ fixes no cusp in characteristic $2$ either. Thus the quotient map is
\'etale, also at $2$, and the genera follow from Riemann--Hurwitz, as $X_1(51)$ and $X_1(63)$
have genus $65$ and $97$.
\end{proof}

\begin{proposition}\label{prop:51_63}
The groups $\torz{51}$ and $\torz{63}$ are not subgroups of the torsion of an elliptic curve
over a sextic field.
\end{proposition}

\begin{proof}
Let $n\in\{51,63\}$ and let $E$ be an elliptic curve over a sextic field $K$ with a point of order
$n$. No elliptic curve over $\F_{2^f}$, $f\leq6$, has a point of order $n$, by the Hasse bound for
$f\leq5$ and because an elliptic curve over $\F_{64}$ cannot have exactly $51$ or $63$ points.
By \Cref{lem:badred} and \Cref{lem:additive}, $E$ has multiplicative reduction at all primes above $2$, so
the trace divisor $x^{(6)}$ reduces to an effective cuspidal divisor modulo $2$, and so does its
push-forward $D$ to $Y=Y_n=X_1(n)/H$, where $H$ is the subgroup generated by $\dm{16}$ for $n=51$ and by $\dm8$ for $n=63$. By
\Cref{lem:quotients_51_63}, $H$ acts freely on $X_1(n)$, also in characteristic $2$, so $Y$
satisfies the assumptions of \Cref{prop:cartier}. The divisor $D$ is an effective
$\Q$-rational divisor of degree $6$ with non-cuspidal support.

The Jacobian $J(Y)$ is isogenous to the part of $J_1(n)$ on which $H$ acts trivially, whose
isogeny factors correspond to the newforms of level dividing $n$ whose characters are trivial on
$H$. For $n=51$ we have $\rk J_1(51)(\Q)=0$ by \cite[Theorem 3.1(2)]{DEHMZ}. For $n=63$ the only
newform of level dividing $63$ with positive analytic rank is \lmfdbnewform{63}{2}{s}{a}, whose
character $\chi$ satisfies $\chi(8)=-1$, as we \gitlink{rank_quotient_63_character.m}{check}, so it
does not contribute to $J(Y_{63})$. Hence in both
cases $J(Y)(\Q)$ is finite, and we apply \Cref{prop:cartier} with $r=1$ and $\psi_1=1$.

For $Y_{51}$ we \gitlink{cusp51_direct_cartier.m}{compute} that
$U=(\Cart-1)H^0(Y_{\F_2},\Omega^1)$ has dimension $30$, that the closed cusps of $Y_{\F_2}$
have \gitlink{cartier_direct_51_union_continuation.m}{degrees} $1^{16},4^4$, so there are
$\binom{21}{6}+4\binom{17}{2}=54808$ effective cuspidal divisors of degree $6$, and that the
restriction of $U$ to each of them has \gitlink{cusp51_retained_rr_kernels_continuation.m}{rank}
$6$. For $Y_{63}$ we
\gitlink{cusp63_direct_cartier.m}{compute} that the closed cusps have degrees $1^9,2^3,3^9,6$, giving $6542$
effective cuspidal divisors of degree $6$, that $U$ has dimension $45$, and that the restriction
of $U$ to each of the $6542$ divisors has rank $6$. By \Cref{prop:cartier} no reduction of $D$ is
possible, which is a contradiction.
\end{proof}

\section{The groups \texorpdfstring{$\torz{57}$}{Z/57Z} and \texorpdfstring{$\torz{65}$}{Z/65Z}}\label{sec:57_65}

For $n=57$ and $n=65$ the Jacobian $J_1(n)$ has positive rank, and over $\F_{64}$ there are
elliptic curves with a point of order $n$, namely supersingular curves with $57$ and $65$ points.
We first show that the corresponding reduction cannot occur.

\begin{proposition}\label{prop:good2}
Let $n\in\{57,65\}$, and let $E$ be an elliptic curve over a sextic field $K$ with a point of
order $n$. Then $E$ has bad reduction at every prime of $K$ above $2$.
\end{proposition}

\begin{proof}
There are no points of order $n$ over number fields of degree at most $3$ \cite{DEHMZ}, so the
point $x\in Y_1(n)(K)$ has degree $6$ and its trace divisor $D$ is a closed point of degree $6$.
Suppose that $E$ has good reduction at a prime $\wp$ of $K$ above $2$. By the Hasse bound $\F_\wp=\F_{64}$,
so $\wp$ is the only prime of $K$ above $2$, and it is unramified. The reduction $\widetilde E$
has a point of order $n$ over $\F_{64}$, so $\#\widetilde E(\F_{64})=n$, and the trace of
Frobenius is $8$ or $0$. Thus $\widetilde E$ is supersingular and $j(\widetilde E)=0$. As
$\widetilde E$ is supersingular, $[2]$ is purely inseparable \cite[Theorem V.3.1]{silverman}, so
$[8]$ and the Frobenius endomorphism $\pi$ of $\widetilde E$ are purely inseparable isogenies of
degree $64$. By \cite[Corollary II.2.12]{silverman} both are the $64$-th power Frobenius
morphism followed by an isomorphism, so $\pi=\zeta\circ[8]$ for an automorphism $\zeta$ of
$\widetilde E$. As $\pi$ acts trivially on
$\widetilde E(\F_{64})$, which has odd order $n$, the automorphism $\zeta$ acts on it as
multiplication by $8^{-1}$. In particular $\zeta$ maps $8\widetilde P$ to $\widetilde P$, where
$\widetilde P$ is the reduction of the point of order $n$, so $(\widetilde E,8\widetilde P)\simeq
(\widetilde E,\widetilde P)$ and $\dm8$ fixes the reduction of $x$. As $\dm8$ is defined over $\Q$, it
also fixes the conjugates of this point, so the reduction of $D$ is fixed by $\dm8$. We \gitlink{check_good2_57_65.m}{check} this also by an
explicit enumeration.

The only newforms of level dividing $n$ with positive analytic rank are
\lmfdbnewform{57}{2}{a}{a} and \lmfdbnewform{65}{2}{a}{a}, which have trivial character, so
$(\dm8-1)J_1(n)(\Q)$ is finite. The class of the degree $0$ divisor $(\dm8-1)D$ is
\[
[(\dm8-1)D]=(\dm8-1)[D-6c_0]+6[\dm8c_0-c_0],
\]
which is torsion, as the first summand lies in $(\dm8-1)J_1(n)(\Q)$ and the second is torsion by
Manin--Drinfeld. It reduces to $0$ modulo $\wp$, as $(\dm8-1)$ kills the reduction of $D$. By
Parent's result $2(\dm8-1)D=\ddiv(f)$ for some $f\in\Q(X_1(n))$, of degree at most $12$, while
Abramovich's bound gives $\gon_\C X_1(n)>14.28$ and $19.99$, respectively. Hence $f$ is constant
and $\dm8D=D$. The least degree of a CM point
on $X_1(n)$ is $12$, respectively $24$, so the support of $D$ is non-CM and, by \Cref{lem:orbits}, $\dm8$ acts freely on its geometric
points.

For $n=57$ the operator $\dm8$ has order $3$ in $\Delta(57)$, so the image of $x$ on
$X_1(57)/\dm8$ has degree $2$, and the image of $x$ on $X_0(57)$ has degree $1$ or $2$. The
curve $X_0(57)$ has no rational non-cuspidal points, and its quadratic points were determined by
Box \cite[\S 4.4]{Box}. The non-CM ones are defined over $\Q(\sqrt{-23})$ or $\Q(\sqrt{13})$, and
the corresponding quadratic field is contained in $K$. The prime $2$ splits in
$\Q(\sqrt{-23})$, which contradicts the uniqueness of the prime of $K$ above $2$. The two
$j$-invariants over $\Q(\sqrt{13})$ are
\[
\frac{3387888351672962316333\sqrt{13}-12215205167504087643323}{2}\quad\text{and}\quad
\frac{-17787\sqrt{13}-61763}{2},
\]
and writing them in the basis $1,t$ of the ring of integers, with
$t=(1+\sqrt{13})/2$, they are congruent to $t$ modulo the prime above $2$, which is inert. As
$t^2+t+1\equiv0$ modulo this prime, these $j$-invariants are non-zero modulo $2$, contradicting
$j(\widetilde E)=0$.

For $n=65$ the operator $\dm8$ has order $2$, so the image of $x$ on $X_1(65)/\dm8$ has degree $3$
and its image on $X_0(65)$ has degree $1$ or $3$. Again $X_0(65)$ has no rational non-cuspidal
points, and by \cite[Theorem 1.1 and \S 5.2.4]{BGG} all its cubic points are defined over the
field $\Q(\alpha)$ with $\alpha^3-\alpha^2+\alpha-2=0$. This field would be contained in $K$, but
it has discriminant $-83$ and the defining polynomial factors modulo $2$ as $T(T^2+T+1)$, so $2$
has two primes above it in $\Q(\alpha)$, which is a contradiction.
\end{proof}

\begin{lemma}\label{lem:57two}
The diamond operator $\dm2$ acts trivially on $J_1(57)(\Q)[2]$.
\end{lemma}

\begin{proof}
Let $J=J_1(57)$. By the Eichler--Shimura relation at $5$ and the injectivity of reduction modulo
$5$ on $J(\Q)[2]$, the group $J(\Q)[2]$ is contained in the kernel of $A_5=T_5-5\dm5-1$ on $J[2]$.
We \gitlink{rational_two_torsion.m}{compute} the action of the Hecke operators on
$H_1(J(\C),\Z)/2H_1(J(\C),\Z)\simeq J[2]$, a space of dimension $170$ over $\F_2$, and find that
the kernel of $A_5$ has dimension $2$ and that $\dm2-1$ vanishes on it. Here we obtain the
diamond operators over $\Z$ from the relation $q\dm q=T_q^2-T_{q^2}$, with $q=2$ and $q=5$.
\end{proof}

\begin{proposition}\label{prop:57}
The group $\torz{57}$ is not a subgroup of the torsion of an elliptic curve over a sextic field.
\end{proposition}

\begin{proof}
Let $E$ be an elliptic curve over a sextic field $K$ with a point of order $57$, and let $D$ be the
trace divisor. By \Cref{prop:good2} and \Cref{lem:additive} the reduction of $D$ modulo $2$ is cuspidal. The
Galois orbits of cusps of $X_1(57)$ over $\Q$ have degrees $1^{18},2^9,18^2$, the orbits of
degree $2$ are defined over $\Q(\zeta_3)$, in which $2$ is inert, and the orbits of degree $18$
have residue degrees $9$ or $18$ modulo $2$. It follows that the reduction of $D$ lifts to a $\Q$-rational
cuspidal divisor $C$ of degree $6$.

Let $g=\dm2$ and $y=(g-1)[D-C]$. As in the proof of \Cref{prop:good2}, $(g-1)J_1(57)(\Q)$ is
finite, so $y$ is a rational torsion class with trivial reduction modulo $2$, and $2y=0$ by Parent's
result. By \Cref{lem:57two} we have $gy=y$. As $g$ has order $9$ in $\Delta(57)$,
\[
0=(g^9-1)[D-C]=(1+g+\cdots+g^8)y=9y=y.
\]
So $(g-1)(D-C)$ is principal, and since its polar divisor has degree at most $12$ and
$\gon_\C X_1(57)>14.28$, it is $0$. Hence $gD=D$, and as the support of $D$ is non-CM,
\Cref{lem:orbits} gives $9\mid\deg D=6$, which is a contradiction.
\end{proof}

For $n=65$ we use \Cref{prop:cartier} on the quotient $Y=X_1(65)/\dm4$. We need an explicit
model of $Y$ over $\F_2$. The quotient of $Y$ by $\dm{51}$ is the normalization of the fibre
product of $X_1(5)$ and $X_1(13)/\dm3$ over the $j$-line. This fibre product only remembers that
the two universal elliptic curves have the same $j$-invariant, and forgets the relative sign,
that is, the isomorphism between them, which $Y$ restores. Over $\F_2$ quadratic twists are
given by Artin--Schreier extensions, so $Y$ is an Artin--Schreier double cover of the fibre
product, which we now make explicit. We use the coordinate $t$ on $X_1(5)$ for which the universal curve
in Tate normal form has $b=c=t$, and a coordinate $z$ on the genus $0$ curve
$X_1(13)/\dm3$. Over $\F_2$ the $j$-maps are
\[
j_5(t)=\frac{t^{12}+t^8+t^4+1}{t^7+t^6+t^5},\qquad
j_{13}(z)=\frac{z^{16}+z^4+1}{z^{28}+z^{25}+z^{24}+z^{21}+z^{20}+z^{17}+z^{16}+z^{13}}.
\]
Over a field of characteristic $2$, an ordinary elliptic curve with $j$-invariant $j\neq0$ is
isomorphic to $y^2+xy=x^3+ax^2+j^{-1}$ for some $a$, which is determined up to adding an element
of the form $w^2+w$. For the universal curve over $X_1(5)$ we have $a_5=t^2/(t+1)^3$, and taking
the trace of the corresponding invariant of the universal curve over $X_1(13)$ under the diamond
operators $\dm3$ gives $a_{13}=(z^7+z^6+z^2+z)/(z^4+z+1)$. We
\gitlink{setup13_twist_parameter_certified.m}{compute} the coordinate $z$, the map $j_{13}$ and
$a_{13}$ from a model of $X_1(13)$ over $\F_2$.

\begin{lemma}\label{lem:Y65}
Let $Y=X_1(65)/\dm4$. The diamond operator $\dm4$ has order $3$ and acts freely on $X_1(65)$, also
in characteristic $2$, so $Y$ has good reduction at $2$ and genus $41$. The function field of
$Y_{\F_2}$ is generated by $t$, $z$ and $w$ with
\[
j_5(t)=j_{13}(z),\qquad w^2+w=a_5+a_{13},
\]
and the residual diamond group $\Delta(65)/\dm4$, of order $8$, is generated by the automorphisms
\[
\tau_T(t,z,w)=\Bigl(\frac1t,\ z,\ w+\frac{t}{t+1}\Bigr),\qquad
\tau_U(t,z,w)=(t,\ z+1,\ w+z),
\]
which have order $4$ and commute. More precisely, $\tau_T$ is the reduction of $\dm3$ or of
$\dm{3\cdot51}$, $\tau_U$ is the reduction of $\dm6$ or of $\dm{6\cdot51}$, and
$\tau_T^2=\tau_U^2$ is the reduction of $\dm{51}$, which acts by $(t,z,w)\mapsto(t,z,w+1)$.
\end{lemma}

\begin{proof}
The element $-4\equiv61\pmod{65}$ has order $3$, and it is $\equiv1\pmod5$. A geometric point
fixed by $\dm4$ would give an automorphism $u\neq\pm1$ of an elliptic curve with $u^3=1$ fixing a
point $P_5$ of order $5$, and then $3P_5=(u^2+u+1)P_5=0$, which is impossible. As
$61\not\equiv\pm1\pmod{13}$, it follows from \Cref{sec:prelim} as in \Cref{lem:quotients_51_63}
that no cusp is fixed either, also in characteristic $2$, so the quotient map is \'etale and
Riemann--Hurwitz gives $2g(Y)-2=(2\cdot121-2)/3$, so $g(Y)=41$. The quotient
$X_1(65)/\dm{4,51}$, where $51\equiv1\pmod5$ and $51\equiv-1\pmod{13}$, is the normalization of
the fibre product of $X_1(5)$ and $X_1(13)/\dm3$ over the $j$-line, as it parametrizes
elliptic curves with a point of order $5$ up to sign and a $\dm3$-orbit of points of order $13$
up to sign. The double cover $Y$ of it recovers the relative sign, which is the choice of an
isomorphism between the two universal curves, and over $\F_2$ the isomorphisms between
$y^2+xy=x^3+a_5x^2+j^{-1}$ and $y^2+xy=x^3+a_{13}x^2+j^{-1}$ correspond to the solutions of
$w^2+w=a_5+a_{13}$. The formulas for $\tau_T$ and $\tau_U$ follow from the identities
$a_5(t^{-1})+a_5(t)=(t/(t+1))^2+t/(t+1)$ and $a_{13}(z+1)+a_{13}(z)=z^2+z$. We
\gitlink{cusp65_genus41_rr_transport_full.m}{verify} that the fibre product is irreducible, that $Y$ has genus $41$ and that
$\tau_T$ and $\tau_U$ are automorphisms of order $4$ that commute.

It remains to identify $\tau_T$ and $\tau_U$ with diamond operators. The maps from
$X_1(65)/\dm{4,51}$ to $X_1(5)$ and to $X_1(13)/\dm3$ commute with the diamond operators, so
$\dm a$ acts on $t$ through $a\bmod5$ and on $z$ through $a\bmod13$. On $X_1(5)$ the operators
$\dm2=\dm3$ act by $t\mapsto-1/t$, hence by $t\mapsto1/t$ over $\F_2$, as we
\gitlink{check_65_diamonds.m}{compute}, and on $X_1(13)/\dm3$ over $\F_2$ the operator $\dm2$
acts by $z\mapsto z+1$, as we \gitlink{setup13_twist_parameter_certified.m}{compute} from the
moduli interpretation. As $3\equiv-2\pmod5$ and $\dm3$ acts trivially on $X_1(13)/\dm3$, while
$6\equiv1\pmod5$ and $\dm6=\dm2\dm3$ on $X_1(13)$, the operators $\dm3$ and $\dm6$ act on
$X_1(65)/\dm{4,51}$ by
$(t,z)\mapsto(1/t,z)$ and $(t,z)\mapsto(t,z+1)$. Two automorphisms of $Y$ over the same
automorphism of $X_1(65)/\dm{4,51}$ differ by an automorphism of $Y$ over $X_1(65)/\dm{4,51}$,
that is, by $1$ or by $\dm{51}$, which acts by $w\mapsto w+1$. Hence $\tau_T$ is $\dm3$ or $\dm3\dm{51}$, and $\tau_U$ is $\dm6$ or
$\dm6\dm{51}$. Finally $\tau_T^2=\tau_U^2=(t,z,w)\mapsto(t,z,w+1)$, and the group generated by
$\tau_T$ and $\tau_U$ has order $8$, so it is the whole residual diamond group.
\end{proof}

\begin{proposition}\label{prop:65}
The group $\torz{65}$ is not a subgroup of the torsion of an elliptic curve over a sextic field.
\end{proposition}

\begin{proof}
Let $E$ be an elliptic curve over a sextic field $K$ with a point of order $65$, and let $D$ be the
push-forward to $Y=X_1(65)/\dm4$ of the trace divisor. By \Cref{prop:good2} and \Cref{lem:additive},
$D$ reduces to an effective cuspidal divisor $\overline C$ of degree $6$ on $Y_{\F_2}$. The closed
cusps of $Y_{\F_2}$ have degrees $1^8,2^4,4^4$, so there are $3576$ effective cuspidal divisors of
degree $6$.

The only newform of level dividing $65$ with positive analytic rank has trivial character, so
$\psi_1=\tau_T-1$ and $\psi_2=\tau_U-1$ have finite image on $J(Y)(\Q)$, where by
\Cref{lem:Y65} we view $\tau_T$ and $\tau_U$ as the diamond operators of $Y$ that they are the
reductions of. We \gitlink{independent/cartier_65.m}{compute}, from the equations of
\Cref{lem:Y65}, that $\psi_1^*H^0(Y_{\F_2},\Omega^1)+\psi_2^*H^0(Y_{\F_2},\Omega^1)$ has dimension
$35$, that its image $U$ under $\Cart-1$ has dimension $28$, and that the restriction of $U$ to
each of the $3576$ divisors has rank $6$. By \Cref{prop:cartier} no such $\overline C$ is the reduction of $D$, which
is a contradiction.
\end{proof}

\section{The sporadic groups and the proof of the main theorem}\label{sec:sporadic}

We now describe the four sporadic groups. Recall that for $b,c$ in a field, the curve
\begin{equation}\label{eq:tate}
E(b,c):\quad y^2+(1-c)xy-by=x^3-bx^2
\end{equation}
in Tate normal form has the point $(0,0)$, and that the curve
\begin{equation}\label{eq:tate2}
E'(b,c):\quad y^2=x^3+cx^2+(1-b-c)bx
\end{equation}
has the point $P=(0,0)$ of order $2$ and the point $Q=(b,b)$.

\begin{proposition}\label{prop:37}
Let $E_0:y^2=x^3-10395x+444150$, which is a model of \lmfdbec{1225}{b}{2}, let $w$ be a root of
\[
f_{37}=x^6-3150x^5+796635x^4-75770100x^3+3111596775x^2-44606598750x-85333003875,
\]
put $K_{37}:=\Q(w)$ and $\delta:=w^3-10395w+444150$, and let
\[
E_{37}:\quad y^2=x^3-10395\delta^2x+444150\delta^3 .
\]
Then $K_{37}$ is the cyclic totally real sextic field \lmfdbnf{6.6.300125.1}, the point
$(\delta w,\delta^2)$ of $E_{37}$ has order $37$ and $E_{37}(K_{37})_\tors\simeq\torz{37}$.
If $K$ is a sextic field and $E/K$ is an elliptic curve with a point of order $37$, then the pair
$(K,E)$ is isomorphic to $(K_{37},E_{37})$. In particular, no group properly containing
$\torz{37}$ is a sextic torsion group.
\end{proposition}

\begin{proof}
Put $g(x)=x^3-10395x+444150$, so that $E_0:y^2=g(x)$ and $\delta=g(w)$. The map
$(x,y)\mapsto(\delta x,\delta^2y)$ is an isomorphism from $E^*:\delta y^2=g(x)$ to $E_{37}$, which
sends $P=(w,1)$ to $(\delta w,\delta^2)$. We \gitlink{verify_sporadic_curves.m}{compute} that
$K_{37}$ is the cyclic totally real field \lmfdbnf{6.6.300125.1}, and we \gitlink{verify_37.m}{compute}
that $P$ has order $37$ and that the irreducible factors over $\Q$ of the $37$-division polynomial
of $E_0$ have degrees $6,6,6,666$.

As $E^*$ is a quadratic twist of $E_0$, the $x$-coordinates $x_r$ of the points $rP$, for
$1\leq r\leq18$, are $18$ distinct roots of this division polynomial. They lie in $K_{37}$, so they
have degree at most $6$, and hence they are exactly the $18$ roots of the three sextic factors. In
particular $\Q(x_r)=K_{37}$ for every $r$. Writing $rP=(x_r,y_r)$, we have $\delta y_r^2=g(x_r)$
with $y_r\neq0$, as $rP$ does not have order $2$.

Now let $(K,E)$ be as in the statement. By \cite[Proposition 1.4]{DKSS} we have
$j(E)=-9317=j(E_0)$, so $E$ is isomorphic to $E_d:dy^2=g(x)$ for some $d\in K^\times$. A point of
order $37$ on $E_d$ has $x$-coordinate $x'\in K$, which is a root of the $37$-division polynomial
of $E_0$ of degree at most $6$, hence a root of one of the sextic factors. So $K=\Q(x')$, and
there is an isomorphism $\tau:K\to K_{37}$ with $\tau(x')=x_r$ for some $r$. As $d$ and $g(x')$
differ by a non-zero square in $K$, the elements $\tau(d)$ and $g(x_r)=\delta y_r^2$ differ by a
square in $K_{37}$, so $E_d^\tau\simeq E^*\simeq E_{37}$ over $K_{37}$. Thus
$(K,E)\simeq(K_{37},E_{37})$.

Finally, we \gitlink{verify_37.m}{compute} that $f_{37}$ has roots modulo $29$ and modulo $41$, and that at the
corresponding primes of degree $1$ of $K_{37}$, which are of good reduction, the reduction of
$E_{37}$ has exactly $37$ points. As reduction is injective on the torsion of order prime to the
residue characteristic, $\#E_{37}(K_{37})_\tors$ divides $37$.
\end{proof}

\begin{proposition}\label{prop:examples}
\begin{enumerate}
\item Let $K_{25}=\Q(a,b)$ with $a^3-a^2+1=0$ and $b^2+(a^2-2a-1)b+a=0$, which is the field
      \lmfdbnf{6.2.28037.1}, put
      \[
      r=\frac{a^2b-ab+b-1}{a(ab-1)},\qquad s=\frac{ab-b+1}{ab},
      \]
      and let $E_{25}=E(rs(r-1),s(r-1))$ as in \eqref{eq:tate}. Then $(0,0)$ has order $25$ and
      $E_{25}(K_{25})_\tors\simeq\torz{25}$.
\item Let $K_{4,12}=\Q(\alpha,i)$ with $\alpha^3-3\alpha-4=0$, which is the field
      \lmfdbnf{6.0.419904.2}, and let $E_{4,12}$ be the base change of
      \lmfdbec{162}{d}{2}: $y^2+xy+y=x^3-x^2+4x-1$ to $K_{4,12}$. Then
      $E_{4,12}(K_{4,12})_\tors\simeq\tg{4}{12}$.
\end{enumerate}
\end{proposition}

\begin{proof}
The point in (1) is the sextic point on $X_1(25)$ found by van Hoeij \cite{vanHoeij}, and the
curve in (2) was found by Gonz\'alez-Jim\'enez and Najman \cite{GN}. We
\gitlink{verify_known_examples.m}{verify} both claims directly in Magma \cite{magma}.
\end{proof}

\begin{proposition}\label{prop:2_24_points}
Let $K_{2,24}=\Q(w)$ with $w^6-2w^4+2w^2-3=0$, which is the field \lmfdbnf{6.2.231852.1}, and let
$E_{2,24}=E'(b,c)$ as in \eqref{eq:tate2} with
\[
b=\tfrac{1}{32}\bigl(3w^5+7w^4-24w^3-15w^2+21w+18\bigr),\qquad
c=\tfrac{1}{16}\bigl(-7w^5+5w^4+24w^3-5w^2-17w-2\bigr).
\]
\begin{enumerate}
\item The point $Q=(b,b)$ has order $24$, and $E_{2,24}(K_{2,24})_\tors\simeq\tg{2}{24}$.
\item The curve $X_1(2,24)$ has exactly $8$ non-cuspidal closed points of degree $6$. They form
      two orbits of length $4$ under $\Delta(24)$, and they all come from $E_{2,24}$.
\item If $K$ is a sextic field and $E/K$ is an elliptic curve with $\tg{2}{24}\subseteq E(K)$,
      then $(K,E)\simeq(K_{2,24},E_{2,24})$.
\end{enumerate}
The curve $E_{2,24}$ has no complex multiplication, its $j$-invariant generates the cubic
subfield of $K_{2,24}$, which has discriminant $-139$, and its conductor has norm $6$.
\end{proposition}

\begin{proof}
Write $K=K_{2,24}$ and $E=E_{2,24}$. We \gitlink{verify_sporadic_curves.m}{compute} (1) and the last
sentence. Moreover $\#\Aut(K)=2$. Indeed, the defining polynomial of $w$ is $h(x^2)$, where
$h(t)=t^3-2t^2+2t-3$ is strictly increasing (its derivative $3t^2-4t+2$ has no real roots) and
$h(0)<0$. So $h$ has exactly one real root, which is positive, and $h(x^2)$ has exactly two real
roots, which differ by sign. Hence, in a real embedding of
$K$, every automorphism of $K$ sends $w$ to $\pm w$, and $w\mapsto-w$ is an automorphism.

We prove (2). The group $E(K)_\tors\simeq\tg{2}{24}$ contains $16$ points of order $24$, and for
each such $Q'$ there are $2$ points of order $2$ not in $\dm{Q'}$, so there are $32$ pairs
$(P',Q')$ generating $E(K)_\tors$, that is, $32$ embeddings $\tg{2}{24}\hookrightarrow E(K)$. Two
of them give the same point of $Y_1(2,24)$ if and only if they differ by an automorphism of $E$.
As $\Aut(E)=\{\pm1\}$, the pairs $(P',Q')$ and $(-P',-Q')=(P',-Q')$ give the same point, so the
$32$ pairs give $16$ distinct $K$-points of $Y_1(2,24)$. Each of them lies on a closed point of degree $6$, as
$\tg{2}{24}\notin\Phi(d)$ for $d\in\{1,2,3\}$, whose residue field is isomorphic to $K$, and such
a closed point has exactly $\#\Aut(K)=2$ points with values in $K$. Hence $E$ gives at least $8$
closed points of degree $6$. On the other hand, by \Cref{prop:2_24}, for every non-cuspidal
closed point $D$ of degree $6$ there are $\delta\in\Delta(24)$ and $i$ with $\delta D\sim L_i$, and
$|L_i|$ contains at most one effective divisor, so there are at most $2\cdot\#\Delta(24)=8$ such
points. Hence there are exactly $8$, they form two orbits of length $4$, and they all come from
$E$. In particular each of $|L_1|$ and $|L_2|$ contains exactly one effective divisor, and the
$16$ points above are all the $K$-points on these $8$ closed points. As a consistency check, we
\gitlink{match_n24_exact_points.m}{verify} that two of these points reduce modulo $5$ to the two survivors
in the proof of \Cref{prop:2_24}.

For (3), let $(K',E')$ be as in (3). It gives a non-cuspidal closed point of degree $6$ on
$X_1(2,24)$, which is one of the $8$ points of (2), so its residue field is isomorphic to $K$.
Choosing such an isomorphism, we obtain a $K$-point on one of these closed points, which is one
of the $16$ points above, all of which come from $E$. As $X_1(2,24)$ is a fine moduli space,
$(K',E')\simeq(K,E)$.
\end{proof}

\begin{proposition}\label{prop:48}
The group $\torz{48}$ is not a subgroup of the torsion of an elliptic curve over a sextic field.
\end{proposition}

\begin{proof}
Suppose $E_1$ is an elliptic curve over a sextic field $K$ with a point $P_1$ of order $48$. By
\Cref{lem:4n}, the curve $E'=E_1/\dm{24P_1}$ has $\tg{2}{24}\subseteq E'(K)$, so by
\Cref{prop:2_24_points}(3) we may assume that $K=K_{2,24}$ and $E'=E_{2,24}$. Let
$\phi:E_1\to E'$ be the quotient isogeny. Its kernel has order $2$, and
$\ker\hat\phi=\phi(E_1[2])$ is generated by a $K$-rational point of order $2$, so $E_1$ is
isomorphic to one of the three curves $2$-isogenous to $E_{2,24}$ over $K$. We
\gitlink{verify_sporadic_curves.m}{compute} that their torsion groups over $K$ are $\torz{24}$,
$\tg{2}{12}$ and $\torz{24}$, which is a contradiction.
\end{proof}

\begin{proof}[Proof of \Cref{thm:main,thm:sporadic}]
By \Cref{prop:37,prop:examples,prop:2_24_points}, together with \cite[Theorem 1.1]{DS}, every group
in $\mathcal B$ is a sextic torsion group. Conversely, let $G$ be a sextic torsion group. If
$G\notin\mathcal B$, then by \Cref{lem:66} $G$ contains a member of $\mathcal G$.
\Cref{table:summary} lists, for each of the $66$ members of $\mathcal G$, where it is shown that
it is not a subgroup of the torsion of any elliptic curve over a sextic field, which is a
contradiction. Hence $\Phi(6)=\mathcal B$, which is \Cref{thm:main}. \Cref{thm:sporadic} follows
from \Cref{prop:37,prop:examples,prop:2_24_points}.
\end{proof}

\begin{remark}
The sextic points on $X_1(37)$ and $X_1(2,24)$ are sporadic, since both curves have only finitely
many points of degree at most $6$: for $X_1(37)$ this follows from $\gon_\Q X_1(37)=18$
\cite[Theorem 1]{DvH} and Frey's theorem \cite{frey}, and for $X_1(2,24)$ it follows from
$\rk J_1(2,24)(\Q)=0$ and $\gon_\Q X_1(2,24)>6$ (see the proof of \Cref{prop:2_24}) by
\cite[Proposition 2.3]{DS}. The curve
$E_{2,24}$ has conductor of norm $6$, which is very small, as was the case for the curve with
torsion $\tg{2}{18}$ over a quintic field in \cite{N5}. The points of degree $6$ that give the four
sporadic groups are listed in \cite{x1mn} as \xpoint{1.25.6.a}, \xpoint{1.37.6.a},
\xpoint{4.12.6.a} and \xpoint{2.24.6.a}.
\end{remark}

\begin{remark}
We do not know whether $E_{25}$ and $E_{4,12}$ are the only elliptic curves over sextic fields
with torsion $\torz{25}$ and $\tg{4}{12}$, respectively. This would require determining all the
sextic points on $X_1(25)$ and all the cubic points on $X_1(4,12)$ over $\Q(i)$, which we did not
attempt, as it is not needed for \Cref{thm:main}: all the groups properly containing these two
groups are ruled out directly.
\end{remark}

\section{The computations}\label{sec:computations}

Our computations were performed in Magma V2.30-1 \cite{magma}, apart from some auxiliary ones in
Python: the list of groups in \Cref{lem:66}, the values of Abramovich's bound and the orders of
diamond operators that we quote, the queries to the LMFDB, some checks of orbits of cusps, and the
combination and replay of data exported from Magma. They were run on the \emph{Mordell} server
at the Department of Mathematics, University of Zagreb, with an AMD EPYC 9175F CPU and $384$ GB of
RAM. We used the package \texttt{mdmagma} \cite{mdmagma} for the moduli-theoretic operations,
Sutherland's optimized models of $X_1(n)$ \cite{sutherland,sutherlandtables}, the
Derickx--Sutherland models of $X_1(m,n)$ \cite{DSmodels}, and Sutherland's package
\cite{sutherlandmagma} for an independent count of points on modular curves over finite fields.
The Hecke operators on closed points of $X_1(m,n)_{\F_p}$ were computed with the routines of
\cite[\S 10]{N5}, extended to $m\geq3$ by carrying both level points through the isogenies, and
these extensions were tested against the full computation of the isogenies on sample points.
The computations took about $11$ hours of CPU time in total, not counting the independent
recomputations described below, of which the one used in the proof of \Cref{prop:65} took about
$10$ hours.

The completeness of the enumeration of closed points of degree at most $6$, which every sieve
computation assumes, was checked independently of the models by counting the isomorphism classes
of elliptic curves with the required level structure over $\F_{p^d}$ for $d\leq 6$, using
$\#Y_1(m,n)(\F_{p^d})=\sum_{e\mid d}e\cdot n_e$, where $n_e$ is the number of non-cuspidal closed points of
degree $e$ (see \githubbare{independent_counts.m}, \githubbare{count_short_sieves.m},
\githubbare{count_42_gl2.m} and \githubbare{verify_35.m} for the curves of \Cref{table:sieve}).

For the check in \Cref{rem:A2_49} we need that $A_2$ annihilates $J_1(49)(\Q)$. By
\cite[Theorem 3.1(2), Corollary 4.14 and Table 2]{DEHMZ} this group is finite with invariant
factors $[7,52367710906884085342]$, so its $2$-primary part is $\{0,z\}$ for an element $z$ of
order $2$. Reduction modulo $2$ is injective on the odd part, and on the reduction of a rational
point Frobenius acts trivially and Verschiebung acts as multiplication by $2$, so $A_2$
annihilates the odd part by the Eichler--Shimura relation. By \Cref{prop:tools}(2) the element
$z$ lies in the kernel of $A_5$ on $J_1(49)[2]$, and we \gitlink{two_torsion_49_T2.m}{compute},
using the action of the Hecke algebra on
$H_1(J_1(49)(\C),\Z)/2H_1(J_1(49)(\C),\Z)\simeq J_1(49)[2]$, that this kernel has dimension $2$
and meets the kernel of $T_2$ trivially. As $T_2$ and $\dm2$ preserve $\{0,z\}$, it follows that
$T_2z=\dm2z=z$, so $A_2z=-2z=0$. We then \gitlink{sieve_49_A2.m}{compute} that $[A_2(D)]\neq0$
for the $8$ representatives of degree $6$, in agreement with the computation with $A_5$.

We checked the computations on which the proofs depend most heavily a second time, with code
written independently of the first implementation. These are the explicit examples of
\Cref{sec:sporadic} (see \githubbare{verify_sporadic_curves.m}), the list of \Cref{lem:66} (see
\githubbare{group_list.py}) and the analytic ranks and character values used in
\Cref{sec:global,sec:formal,sec:57_65}, read directly from the LMFDB (see
\githubbare{lmfdb_ranks.py}). In the folder \githubdir{independent} of the repository we recomputed
the $2$-torsion computations of \Cref{rem:A2_49} and \Cref{lem:57two} (in both conventions for
the action of the Hecke algebra on homology), and, with models of the curves constructed from
scratch, the doubled cuspidal computations of \Cref{prop:doubled_cases} for $X_1(45)$ and
$X_1(3,21)$ and the Cartier operator computations of \Cref{prop:51_63}. For \Cref{prop:65} the
computation in this folder is the one used in the proof, and the original one
(\githubbare{cusp65_genus41_rr_transport_full.m}) gives the same dimensions $35$ and $28$, and
rank $6$ on the $1792$ of these divisors that are compatible with the fields of definition of the
cusps of $X_1(65)$. In all cases the results agree.

\bibliographystyle{amsalpha}
\bibliography{sextic}

\end{document}